\documentclass[11pt]{article}%
\usepackage{amsfonts}
\usepackage{amssymb,amsmath,amsthm, amsfonts}
\usepackage{mathtools}
\usepackage{color}

\makeatletter

\newcommand{\Rmnum}[1]{\expandafter\@slowromancap\romannumeral #1@}
\makeatother

\providecommand{\U}[1]{\protect \rule{.1in}{.1in}}
\newtheorem{theorem}{Theorem}[section]

\newtheorem{corollary}{Corollary}[section]

\newtheorem{definition}{Definition}[section]

\newtheorem{lemma}{Lemma}[section]

\newtheorem{proposition}{Proposition}[section]
\newtheorem{remark}{Remark}[section]
\newtheorem{example}{Example}[section]

\def\sup{\mathop{\rm sup}}
\makeatletter
   
   \@addtoreset{equation}{section}
\makeatother

\begin{document}

\title{\LARGE\bf Mean Field Games of Major-Minor Agents with Recursive Functionals}

\author{
Jianhui Huang \thanks{ Department of Applied Mathematics, The Hong Kong Polytechnic University,
 Hong Kong, China.  Email: {james.huang@polyu.edu.hk}. This author's research is partially supported by RGC Grant PolyU 15301119, 15307621, N PolyU504/19, NSFC 12171407 and KKZT.}
 \and
Wenqiang Li\thanks{Corresponding author. School of Mathematics and Statistics, Shandong University, Weihai 264209, P.R. China.  Email: {wenqiangli@sdu.edu.cn}. This author's research is supported by the NSF of P.R. China (No. 12101537, 12271304).}
\and
Harry Zheng \thanks{Department of Mathematics, Imperial College,
 London SW72BZ, UK.  Email:{h.zheng@imperial.ac.uk}. This author is partially supported by  EPSRC (UK) grant (EP/V008331/1).}
}

\date{\today}
\maketitle

\bigskip
\noindent \textbf{Abstract.} This paper investigates a novel class of mean field games involving a major agent and numerous minor agents, where the agents' functionals are recursive with nonlinear backward stochastic differential equation (BSDE) representations. We term these games ``recursive major-minor" (RMM) problems. Our RMM modeling is quite general, as it employs empirical (state, control) averages to define the weak couplings in both the functionals and dynamics of all agents, regardless of their status as major or minor. We construct an auxiliary limiting problem of the RMM by a novel unified structural scheme combining a bilateral perturbation with a mixed hierarchical recomposition. This scheme has its own merits as it can be applied to analyze more complex coupling structures than those in the current RMM. Subsequently, we derive the corresponding consistency condition and explore asymptotic RMM equilibria. Additionally, we examine the RMM problem in specific linear-quadratic settings for illustrative purposes.

\medskip


\medskip

\noindent \textbf{MSC2020 subject classifications.} 93E20, 60H10, 60K35.

\section{Introduction} Mean field game (MFG) theory was independently introduced by \cite{LL2007} and \cite{HMC2006} from different perspectives, serving as an effective methodology for analyzing controlled large population (LP) systems. Typically, a LP system comprises a large number of agents with interactions through their empirical distribution or averages. These interactions are in a weak sense, as the degrees of couplings among agents diminish rapidly when the number of agents tends to infinity. 
A core element of the MFG theory is the construction of a limiting auxiliary problem, under which all agents can be largely disentangled, allowing for characterizations of some decentralized approximate equilibrium through a consistency matching. As such, MFG analysis can significantly reduce the dimension of controlled LP systems that each agent needs to analyze, greatly simplifying related numerical analysis. A substantial body of research has been dedicated to the MFG theory, yielding fruitful outcomes. A partial list of literature relevant to this current work includes {\cite{BCY2016, BCR2024, CCP2020, CZ2016, HWW2016, H2010, MH2020, NC2013}.} This paper focuses on a new class of MFG problems with the following weakly-coupled LP system including a major $\mathcal{A}_0$ and multiple minors $\{\mathcal{A}_i\}_{i=1}^{N}$, whose states $X^0$ and $\{X^i\}_{i=1}^{N}$ satisfy the following \emph{forward} stochastic differential equations (SDEs)
\begin{equation}\label{major-state}
\left\{
\begin{aligned}
dX_t^0=&b^0(t,X_t^0,{u_t^0}; X_t^{(N)},{u_t^{(N)}})dt+{\sigma^{0}(t,X_t^0,{u_t^0}; X_t^{(N)},{u_t^{(N)}})}dW^0_t,\\
dX_t^i=&{b(t,X_t^i,u_t^i;
X_t^0,u_t^0;X_t^{(N)},u_t^{(N)})dt
+\sigma(t,X_t^i,u_t^i;X_t^0,u_t^0;X_t^{(N)},u_t^{(N)})dW_t^i},\end{aligned}
\right.
\end{equation}with initial conditions $X_0^0=x_0^0\in\mathbb{R}^n, X_0^i=x_0\in\mathbb{R}^n.$ Here, $X_\cdot^{(N)}=\frac1N\sum_{i=1}^NX_\cdot^i,\ u_\cdot^{(N)}=\frac1N\sum_{i=1}^Nu_\cdot^i$
stand for the (forward) state-average and control-average of all minors, respectively, and $W=(W^0,W^1,$ $\cdots,W^N)^\top$ is a $(N+1)$-dimensional standard Brownian motion on a complete probability space $(\Omega,\mathcal{F},{\mathbb{P}})$, where
$W^i$ is the idiosyncratic noise for $\mathcal{A}_{i}$ while $W^0$ is the common noise.   $\{\mathcal{A}_{k}\}_{k=0}^{N}$ aim to maximize recursive-type functionals  $\{J_k\}_{k=0}^N$ given by
\begin{equation}\label{091601}
\begin{aligned}
J_0(u^0; u^1,\cdots,u^N)=&\Gamma^0(Y_0^0)+\mathbb{E}\int_0^Tg^0(t,\Theta_t^0, u_t^0;\Theta_t^{(N)},u_t^{(N)})dt,\\
J_{i}(u^{i}; u^{-i}, u^0)
=&\Gamma(Y_0^i)+\mathbb{E}\int_0^Tg(t,\Theta_t^i,u_t^i;\Theta_t^0, u_t^0;\Theta_t^{(N)},u_t^{(N)})dt,\ i=1,\cdots,N,
\end{aligned}
\end{equation}
where $u^{-i}=\{ \cdots, u^{i-1}, u^{i+1}, \cdots \}$ is the control profile except that of $\mathcal{A}_{i}$; $\Theta_\cdot^{(N)}=(X_\cdot^{(N)},Y_\cdot^{(N)},Z_\cdot^{(N)})$ is state-average triple with $
Y^{(N)}_\cdot=\frac1N\sum_{i=1}^NY_{\cdot}^i$ as the \emph{recursive state} average and $Z^{(N)}_{\cdot}=\frac1N\sum_{i=1}^NZ_{\cdot}^i$ the \emph{intensity state} average (see \cite{CE2002}); the state-triple $\{\Theta_\cdot^{k}=(X_\cdot^{k},Y_\cdot^{k}, {Z}_\cdot^{k})\}_{k=0}^N$ satisfy  \eqref{major-state} and the  \emph{backward} SDE (BSDE) motivated by recursive utilities in economics \cite{DE1992, L2004, LQ2003}:
\begin{equation}\label{major-utility}
\left\{
\begin{aligned}
-dY_t^0=&f^0(t,\Theta_t^0,u_t^0;  \Theta_t^{(N)},u_t^{(N)})dt-Z_t^{0,0}dW_t^0
-\sum_{j=1}^{N}Z_t^{0,j}dW_t^j,\\
-dY_t^i=&f(t,\Theta_t^i,u_t^i;\Theta_t^0,u_t^0; \Theta_t^{(N)},u_t^{(N)})dt
-Z_t^{i,0}dW^0_t-Z_t^{i,i}dW^i_t
-\sum_{j=1, j \neq i}^{N}Z_t^{i,j}dW_t^j,\\
Y_T^0=&\Phi^0(X_T^0,X_T^{(N)})+\xi^0,\ Y_T^i=\Phi(X_T^i,X_T^0,X_T^{(N)})+\xi^i,
\end{aligned}
\right.
\end{equation}
where for $i=1, \cdots, N$, $Z^{i}=(Z^{i,0}, Z^{i,i}, {\{Z^{i,j}\}_{j\neq i,j=1}^N})$ and $(Z^{i,0},Z^{i,i})$ is the  \emph{principal} intensity component, and $\{Z^{i,j}\}_{j \neq i}$ the  \emph{marginal} components. 
Since remainder terms $\sum_{j=1}^{N}Z_t^{0,j}dW_t^j$ and $\sum_{j=1, j \neq i}^{N}Z_t^{i,j}dW_t^j$ vanish as $N \rightarrow \infty$ ({see Remark \ref{re-remainder}}), we focus on the principal terms and set $Z^{(N)}_{\cdot}=\frac1N\sum_{i=1}^N (Z^{i,0}_\cdot,Z^{i,i}_\cdot)$ as averages on principal intensities.

We refer \eqref{major-state}-\eqref{major-utility} as the recursive major-minor (RMM) problem. We defer its detailed assumptions later, and highlight its modeling features first as follows.

\textbf{Modeling features}. (i) The RMM model delves into the interaction between a major $\mathcal{A}_{0}$ and a large number of minors $\{\mathcal{A}_{i}\}_{i=1}^{N}$. Traditional MFG studies assume that agents are all ``minor" or ``negligible", meaning that an individual agent's action cannot significantly impact the behaviors of the population at a macro scale. Associated MFGs are thus referred to as ``symmetric" because it suffices to examine a representative agent, provided agents are homogenous hence statistically exchangeable. In contrast, our RMM explores asymmetric interactions where agents having varying decisional capacities. A major agent may significantly influence the population's behavior through her own decisions, whereas numerous minor agents can only affect the population through collective actions. This model is more realistic than the homogeneous minor setting, as it captures a range of diversified interaction mechanisms, see \cite{BCY2016,CZ2016,H2010,MH2020,NC2013}.

(ii) The RMM model further posits that objectives of all agents are represented recursively through nonlinear BSDEs, such as \eqref{major-utility}, with \emph{non-additive} drivers $f^{0}$ or $f$. The inclusion of recursive functionals in MFG studies is motivated by their advantageous decision-theoretic properties, especially in the current LP context featuring complex decision couplings. Indeed, recursive functionals are well-suited for decision theory due to their capability to explain various observed non-standard decision behaviors, such as the separation of inter-temporal substitution and risk aversion. Consistently, recursive functionals extend classical expected functionals (see \cite{AH2021,BLP2014,EPQ2001}), which are relevant to a special class of BSDEs with \emph{additive} drivers.

(iii) The RMM problem restricts its weak coupling to empirical averages, refraining from discussing a more extensive empirical measures or distributions. Despite this limitation, the weak coupling of RMM remains quite general, as it is integrated into both the dynamics and payoff functionals of the major $\mathcal{A}_{0}$ and all minors $\{\mathcal{A}_{i}\}_{i=1}^{N}$, encompassing elements from both the state and the control. Moreover, due to recursive functionals, the state averages is enriched by including not only $X^{(N)}$ on the objective (forward) states; but also $(Y^{(N)}, Z^{(N)})$ on the (recursive and intensity) states reflecting the subjective averaged-out beliefs. Specifically, the intensity coupling $Z^{(N)}$ characterizes some average on risk (ambiguity) aversion across all agents.

\textbf{Literature comparison.} \cite{CW2017,H2010} and \cite{MH2020} introduced the major-minor MFG within a linear-quadratic-Gaussian (LQG) framework on finite and infinite horizon, respectively. They employed augmented Riccati equations to characterize consistency conditions. \cite{NC2013} extended these major-minor studies to a nonlinear setting using the stochastic Hamilton-Jacobi-Bellman (HJB) approach, where the weak coupling is restrictive; for instance, the major's state cannot enter the dynamics of the minors. Besides, \cite{BCY2016} investigated a class of major-minor MFG problems also through the stochastic HJB approach. \cite{CCP2020} studied major-minor MFG by master equations where the agents take closed loop control. Recently, \cite{BCR2024} investigated a type of MFG problems with asymmetric information between major and minor agents. \cite{CZ2016} explored nonlinear major-minor MFG with general weak couplings, allowing the major's state to enter the dynamics of the minor agents. Additionally, the limiting control problem of the major agent incorporates an endogenous mean-field term, based on an approximation through a two-agent non-zero-sum game. A forward type of maximum principle was utilized in this context.

Our paper distinguishes itself from aforementioned works by its focus on a nonlinear major-minor interaction with recursive functionals, and associated methodology of a backward-forward type of stochastic maximum principle. Our RMM modeling is particularly noteworthy for its introduction and detailed analysis of the weak couplings of the backward (recursive and intensity) state-average $(Y^{(N)}, Z^{(N)})$ originating from the recursive functional. As previously mentioned, these couplings hold significant decision-making impacts and, to the best of our knowledge, have not been systematically addressed in the MFG literature. Consequently, the maximum principle we adopt and the consistency condition we derive take unique forms that differ from those in \cite{CZ2016}. Additionally, we apply our general nonlinear outcomes to specific RMM problems in LQG settings. Our LQG-RMM studies not only recap and extend existing results on the \emph{forward} MFG studies, but also provide new insights into its \emph{backward} counterpart.

Another relevant work is \cite{BLP2014}, which also explored the major-minor interaction and recursive functionals. However, it is framed on a weak formulation, substantially different from our strong formulation. Specifically, the approach of \cite{BLP2014} originates from a variant of the Girsanov transformation and optimization of a Hamiltonian function, whereas our analysis is rooted in a refined backward-forward stochastic maximum principle. More significantly, all minors in \cite{BLP2014} are cooperative, and the associated MFG thus encompasses a two-layer mixed structure: all (cooperative) minor agents form a mean-field team (rather than a game) in an inner layer, while in an outer layer, the interacting major agent and a representative minor agent induce a two-person, non-mean-field game. Although also termed a MFG, the model in \cite{BLP2014} is essentially a hybrid of a mean-field \emph{team} and a non-mean-field-type two-person game, which contrasts with the decision structure we investigate. Consequently, the consistency matching, a central step in MFG analysis, is not applicable in \cite{BLP2014} at all. In contrast, our work is distinguished by a novel scheme to auxiliary control construction and consistency matching, as detailed below. Furthermore, the limiting equilibria in \cite{BLP2014} are characterized as saddle points, which are remarkably different from our non-zero-sum setup. Additionally, the admissible controls in \cite{BLP2014} take feedback forms and are compact, unlike our open-loop and unbounded convex admissibility.

\textbf{A unified structural scheme.} Last but not least, as the core element in MFG, the auxiliary problem of the RMM is formulated through a novel \emph{structural scheme}, which incorporates a bilateral perturbation and a hierarchical recomposition. This scheme facilitates a more incisive auxiliary construction through a sequential network, enabling a clear-cut realization of a complex mixture involving two minor agents: one is exogenous and the other is endogenous, alongside the endogenous major agent. This scheme is new in the MFG literature and distinguishes our work from previous studies, where auxiliary constructions are based on heuristic arguments. More importantly, such scheme offers a unified methodology to tackle more complex LP interactions for which heuristic arguments are no longer tractable. For instance, when the LP system consists of heterogenous agents with varying beliefs on model uncertainties.

\textbf{Contributions.} (i) We introduce a new class of RMM problems featuring major-minor asymmetric interactions and recursive objectives. (ii)  We present a novel unified structural scheme to construct its pivotal auxiliary problem. (iii) We derive a new class of mean-field type of forward-backward SDEs (FBSDEs) to characterize the consistency condition of the RMM problem. (iv) We examine LQG studies in the RMM context in detail to gain deeper insights.

The remainder of this paper is organized as follows: Section 2 introduces basic assumptions of the RMM problem. Section 3 presents a unified structural scheme for the RMM problem, including a bilateral perturbation and a mixed triple-agent two-layer analysis. Section 4 studies the auxiliary control and associated consistency condition (CC) of the RMM problem, and verifies its approximate Nash equilibrium. Section 5 devotes to some LQG-RMM problems. Section 6 concludes, and some technical proofs and heavy notations are given in Appendix.

\section{Preliminary} For $i=0,1,\cdots, N$, let $\mathbb{F}^{i}=\{{\mathcal F}_t^{i}\}$ be the complete filtration generated by the Brownian motion $W^{i}$; namely, ${\mathcal F}_t^{i}=\sigma(W_i(s), 0\leq s\leq t) \vee \mathcal{N}_\mathbb{P}$ with $\mathcal{N}_\mathbb{P}$ the set of all $\mathbb{P}$-null sets in $\mathcal F$. Then, $\mathbb{F}=\{{\mathcal F}_t\}_{0 \leq t \leq T}=\{\sigma(\bigcup_{i=0}^N {\mathcal F}_t^{i})\}_{0 \leq t \leq T}$ denotes the \emph{centralized} information generated by the Brownian motion $W$, and by $\mathbb{F}^{i,0}=\{{\mathcal F}^{i,0}_t\}_{0 \leq t \leq T}=\{\sigma( {\mathcal F}_t^{i}\bigcup {\mathcal F}_t^{0} )\}_{0 \leq t \leq T}$ the \emph{decentralized} information for a generic minor agent $\mathcal{A}_i$, $i=1,\cdots, N$.
Let $U_0\subset\mathbb{R}^{k_0}$ and $U\subset\mathbb{R}^{k}$ be two convex sets.

\begin{definition}
$u^0$ is a centralized (resp. decentralized) admissible control for  $\mathcal{A}_0$, if it is an  $\mathbb{F}$-adapted (resp. $\mathbb{F}^0$-adapted) $U_0$-valued process with ${||u||_{L^{2}}:=}\mathbb{E}\int_0^T|u_t|^2dt<\infty$.
Similarly, for $i=1,2,\cdots,N$, a $U$-valued process $u^i$ is called a centralized (resp. decentralized) admissible control for  $\mathcal{A}_i$, if it is $\mathbb{F}$-adapted (resp. $\mathbb{F}^{i,0}$-adapted) with $||u^{i}||_{L^{2}}<\infty$.
Let {$\mathcal{U}_0^c$ and $\mathcal{U}^c$} (resp. $\mathcal{U}_0^d$ and $\mathcal{U}_i^d$) be the set of all centralized (resp. decentralized) controls for $\mathcal{A}_0$ and $\mathcal{A}_i$, respectively.

%
\end{definition}
\begin{definition}\label{def1}
For any $\varepsilon>0$, we say a $(N+1)$-tuple admissible controls $(u^{0,*},u^{1,*},\cdots,u^{N,*})
\in\mathcal{U}^c_0\times\mathcal{U}^c\times\cdots\mathcal{U}^c$ (resp. $\in\mathcal{U}^d_0\times\mathcal{U}^d_1\times\cdots\mathcal{U}^d_N$) depending on $\varepsilon$ is a centralized (resp. decentralized) approximate $\varepsilon$-Nash equilibrium, if for all $(u^{0},u^{1},\cdots,u^{N})
\in\mathcal{U}^c_0\times\mathcal{U}^c\times\cdots\mathcal{U}^c$ (resp. $\in\mathcal{U}^d_0\times\mathcal{U}^d_1\times\cdots\mathcal{U}^d_N$), we have
\begin{equation}\nonumber
\begin{aligned}
&J_0(u^{0,*},u^{1,*},\cdots,u^{N,*})\geq J_0(u^{0},u^{1,*},\cdots,u^{N,*})-\varepsilon,\\
&J_i(u^{0,*},u^{1,*},\cdots,u^{N,*})\geq J_i(u^{0,*},\cdots,u^{i-1,*},u^{i},u^{i+1,*},\cdots,u^{N,*})-\varepsilon,\ i=1,\cdots,N.
\end{aligned}
\end{equation}\end{definition}
The exact Nash equilibrium corresponds to the case when $\varepsilon=0$. Now we impose some assumptions on the following coefficients of \eqref{major-state}-\eqref{major-utility} of the RMM problem:\begin{equation}\nonumber
\begin{aligned}
&(b^0,\sigma^0): \ \ [0,T] \times \mathbb{R}^{n}\times
U_0\times
 \mathbb{R}^{n}\times U \rightarrow\mathbb{R}^n, \quad (b,\sigma): \ \ [0,T] \times\mathbb{R}^{n}\times
U\times \mathbb{R}^{n}\times
U_0\times
 \mathbb{R}^{n}\times U \rightarrow\mathbb{R}^n,\\
&(f^0,g^0): \ \ [0,T]\times \mathbb{R}^{n+m+m} \times U_0\times
 \mathbb{R}^{n+m+2m}\times
U\rightarrow\mathbb{R}^m\times\mathbb{R},\\
  &(f,g): \ \ [0,T]\times
 \mathbb{R}^{n+m+2m} \times
U\times \mathbb{R}^{n+m+m}\times U_0\times
 \mathbb{R}^{n+m+2m}\times
U\rightarrow\mathbb{R}^m\times\mathbb{R},\\
 &\Phi^0:\  \mathbb{R}^{n+n}\rightarrow{\mathbb{R}^m},\quad \ \ \Phi:\  \mathbb{R}^{n+n+n}\rightarrow{\mathbb{R}^m},\quad \ \
 {(\Gamma^0,\Gamma):\  \mathbb{R}^m\rightarrow\mathbb{R}\times\mathbb{R}},
\end{aligned}
\end{equation}
\textbf{{Assumption (A1)}}
(i)   \emph{$(b^0,\sigma^0)$ and $(b,\sigma)$ are continuously differentiable in $(x^0,u^0;x^{(N)},u^{(N)})$ and $(x,u;x^0,u^0;$ $x^{(N)},u^{(N)})$, respectively. All the derivatives of $(b^0,\sigma^0,b,\sigma)$ are bounded.}\\
(ii) \emph{ $(f^0,g^0)$ and $(f,g)$ are continuously differentiable in $(\theta^0,\theta^{(N)})$ and $(\theta,\theta^0,\theta^{(N)})$, respectively, where
$
\theta^0=(x^0,y^0,z^0,u^0),\ \theta^{(N)}=(x^{(N)},y^{(N)},z^{(N)},u^{(N)}),\ \theta=(x,y,z,u).$
$\Phi^0$, $\Phi$ and $(\Gamma^0,\Gamma)$ are continuously differentiable in $(x^0,x^{(N)})$, $(x,x^0,x^{(N)})$ and $y$, respectively.}\\
 (iii) \emph{All derivatives of $(f^0,f,\Phi^0,\Phi)$ are bounded.
 The derivatives of
$g^0,g$ and $(\Gamma^0,\Gamma)$ are bounded by $C(1+|\theta^0|+|\theta^{(N)}|)$,  $C(1+|\theta|+|\theta^0|+|\theta^{(N)}|)$ and $C(1+|y|)$, respectively, for some $C>0$.} \\
(iv) \emph{ $\xi^0$ and $\xi^i$, $1\leq i\leq N$, are square-integrable, $\mathcal{F}_T^{0}$- and $\mathcal{F}_T^{i,0}$-measurable, respectively. Moreover,
 $\xi^1,\cdots,\xi^N$ are independent and identically distributed conditionally on $\mathcal{F}_T^0$.}\\
(v) \emph{ $(b^0,\sigma^0,b,\sigma;f^0,f,g^0,g)$ are uniformly continuous in $t$.
}

For each admissible $(N+1)$-tuple $\{u^j\}_{j=0}^N$, the coupled SDE system \eqref{major-state} (resp. BSDE \eqref{major-utility}) admits a unique solution tuple $\{X^j\}_{j=0}^N$ (resp. $\{Y^j,Z^j\}_{j=0}^N$) under  (A1), and the recursive-type payoffs $\{J_i\}_{i=0}^N$ are well-defined.

\section{A unified structural scheme of the RMM problem}
For a fixed $N$, if each agent has access to centralized information about all agents, including their instantaneous states realized and controls adopted, the RMM becomes a classical but high-dimensional $(N+1)$-agent game. Existence or uniqueness of its \emph{exact} Nash equilibrium(s) can be ensured under certain mild but high-dimensional conditions, including semi-continuity, coercivity and concavities on $ \{J_i\}_{0=1}^{N}$, or compactness on control admissibility. The open-loop (exact) equilibrium, denoted as $\{u^{i,*}\}_{i=0}^{N} \in {\mathcal{U}_{0}^{c}}\times\prod_{i=1}^N\mathcal{U}_i^c$, can be further characterized through a system of $N+1$ stationary conditions. However, this procedure to exact equilibriums is only feasible in theory and becomes impractical due to the curse of dimensionality when $N$ is large.

MFG theory offers one resolution to constructing near-optimal decentralized strategies, as an alternative approximation to the exact Nash equilibria. A key challenge in MFG is the construction of an auxiliary problem for dimension reduction. Previous MFG studies have constructed auxiliary problems intuitively based on heuristic arguments, effective only when the underlying coupling structure is not overly complex. For instance, when all agents are symmetric as in {\cite{DHW2018, HWW2016, HMC2006, LL2007}}, or even an asymmetric dominant major agent is included as in \cite{BCY2016, BCR2024, CCP2020, CW2017, CZ2016, H2010, MH2020, NC2013}. However, heuristic analysis becomes infeasible to analyze LP systems with more intricate couplings. One reason is that it fails to effectively configure a complex logic network in which various representative agents shall be mutually connected by an ``exogenous---endogenous" relation. Alternatively, we propose a structural scheme that can not only well amount for the extreme generality of the weak coupling in RMM problems, but also lay down an unified foundation to analyze more general and complex LP couplings. In current RMM context, this scheme yields a bilateral perturbation and a triple-agent two-layer game, as discussed below.

\subsection{A bilateral perturbation: the major agent} Letting $\varepsilon=0$ in Definition \ref{def1}, $\mathcal{A}_0$ faces an optimization problem: $\sup_{u^0\in\mathcal{U}_0^c} J_{0}(u^0,u^{1,*},\cdots,u^{N,*}),$ by assuming that all minor $\{\mathcal{A}_i\}_{i=1}^{N}$ implement \emph{exact} Nash equilibrium {$\{u^{i,*}\}_{i=1}^{N}$. When $\mathcal{A}_0$ adopts a perturbed \emph{centralized} control $u^0\in\mathcal{U}_0^c$ instead the exact $u^{0,*}$, her state \eqref{major-state} becomes
\begin{equation}\label{N-major-state}
\left\{
\begin{aligned}
d{X}_t^{0,\dagger}=&b^0(t,{X}_t^{0,\dagger},u_t^0; X_t^{(N),\dagger},u_t^{(N),\ast})dt+\sigma^0(t,{X}_t^{0,\dagger},u_t^0; X_t^{(N),\dagger},u_t^{(N),\ast})dW_t^0,\\
dX_t^{i,\dagger}=&b(t,X_t^{i,\dagger},u_t^{i,*};{X}_t^{0,\dagger},u_t^0; X_t^{(N),\dagger},
u_t^{(N),*})dt+\sigma(t,X_t^{i,\dagger},u_t^{i,*};{X}_t^{0,\dagger},u_t^0; X_t^{(N),\dagger},
u_t^{(N),*})dW_t^i,
\end{aligned}
\right.
\end{equation}
with ${X}_0^{0,\dagger}=x_0^0,\ X_0^{i,\dagger}=x_0,\ i=1,\cdots,N.$ Here, $X_\cdot^{(N),\dagger}:=\frac1N\sum_{i=1}^NX_\cdot^{i,\dagger}$
%
is a \emph{quasi-realized} state average with the superscript ``$^\dag$" to emphasize its dependence on the major's perturbed $u^0$; whereas $u_\cdot^{(N),\ast}:=\frac1N\sum_{i=1}^Nu_\cdot^{i,*}$ is the \emph{exact-realized} control average only depending on the exact strategies so the superscript ``$^\ast$" is still applied. This is essentially an open-loop feature.

By ``quasi-realized," the states are not ``exactly" the ones to be realized when all agents apply their exact strategies; instead, they are ``quasi-exact" as only $\mathcal{A}_0$ deviates from the exact one by adopting a \emph{perturbed} control. Actually, an exact Nash equilibrium $\{u_t^{i,*}\}_{i=1}^{N}$, in its open-loop sense, are defined directly on the basic inputs {$(t, \omega)$} rather than on the ``intermediate" states. Thus, a perturbed $u^0$ by $\mathcal{A}_0$, will not change the controls of the minors. This is very different to the closed-loop case, for which a perturbed $u^0$ will change the major's state, so further alter the implementation of $u^{i,*}$ constructed on these realized states.
Using the notation
\begin{equation*}
\begin{aligned}
&Y_{t}^{(N),\dagger}=\frac1N\sum_{i=1}^NY_{t}^{i,\dagger},\quad Z_t^{(N),\dagger}=\frac1N\sum_{i=1}^NZ_t^{i,\dagger}
=\frac1N\sum_{i=1}^N(Z_t^{i,\dagger,0},
Z_t^{i,\dagger,i},\{Z_t^{i,\dagger,j}\}_{j\neq i})^\top,\\
&\Theta_t^{j,\dagger}=(X_t^{j,\dagger},Y_t^{j,\dagger},Z_t^{j,\dagger}),\ j=0,1,\cdots,N,\quad \Theta_t^{(N),\dagger}=(X_t^{(N),\dagger},Y_t^{(N),\dagger},Z_t^{(N),\dagger}),
\end{aligned}
\end{equation*}and similar to \eqref{N-major-state},
we can get the following quasi-realized coupled BSDEs
  \begin{equation}\label{fu-3}
\left\{
\begin{aligned}
-dY_t^{0,\dagger}=&{f^0\big(t,u_t^0,\Theta_t^{0,\dagger};  u_t^{(N),\ast},\Theta_t^{(N),\dagger}\big)}dt
-Z_t^{0,\dagger}dW_t^0-\sum_{j=1}^{N}Z_t^{0,\dagger,j}dW_t^j,\\
-dY_t^{i,\dagger}=&f\big(t,{\Theta_t^{i,\dagger},u_t^{i,*};
\Theta_t^{0,\dagger},u_t^0;
\Theta_t^{(N),\dagger},u_t^{(N),*}}\big)dt
{-Z_t^{i,\dagger,0}dW_t^0-Z_t^{i,\dagger,i}dW_t^i}
-\sum_{j=1, j \neq i}^{N}Z_t^{i,\dagger,j}dW_t^j,\\
Y_T^{0,\dagger}=&\Phi^0\big(X_T^{0,\dagger},X_T^{(N),\dagger}\big)+\xi^0,\quad Y_T^{i,\dagger}=\Phi\big(X_T^{i,\dagger};X_T^{0,\dagger},X_T^{(N),\dagger}\big)
+\xi^i,\quad
i=1,\cdots,N.\\
\end{aligned}
\right.
\end{equation}
We aim to analyze the asymptotic limit as $N\rightarrow\infty$ with $u_t^{(N),*}\rightarrow\overline{u}_t^*\in\mathcal{F}_t^0.$
Then we take the limit of \eqref{N-major-state}, \eqref{fu-3} and the related recursive payoff for the major agent.
By the continuity of coefficients $b^0$ and $b$,
$
\lim_{N\rightarrow\infty}
b^0(t,{X}_t^{0,\dagger},u_t^0; X_t^{(N),\dagger},u_t^{(N),\ast})
=
b^0(t,\mathbb{X}_t^{0,\dagger},u_t^{0};
\overline{\alpha}_t^{\dagger},
\overline{u}_t^{*}),
$ and
\begin{equation}\nonumber
\begin{aligned}
&\lim_{N\rightarrow\infty}\frac1N\sum_{i=1}^N
b(t,X_t^{i,\dagger},u_t^{i,*};{X}_t^{0,\dagger},u_t^0; X_t^{(N),\dagger},
u_t^{(N),*})\\
=&\lim_{N\rightarrow\infty}\frac1N\sum_{i=1}^N
\left[{\lim_{N\rightarrow\infty}\frac1N\sum_{i=1}^N}b\left(t,a;
{X}_t^{0,\dagger},
u_t^0;X_t^{(N),\dagger}, u_t^{(N),*}\right)\right]_{ a=(X_t^{i,\dagger}, u_t^{i,*})}\\
\end{aligned}
\end{equation}
\begin{equation}\nonumber
\begin{aligned}
=&\lim_{N\rightarrow\infty}\frac1N\sum_{i=1}^N
b(t,X_t^{i,\dagger}, u_t^{i,*};{\mathbb{X}}_t^{0,\dagger},
u_t^0;\overline{\alpha}_t^{\dag},
\overline{u}_t^{*})
={\mathbb{E}_{t}\Big[b(t,{\mathbb{X}}_t^{1,\dagger},u_t^{1,*};
\mathbb{X}_t^{0,\dagger},u_t^0;\overline{\alpha}_t^{\dag},
\overline{u}_t^{*})\Big]},
\end{aligned}
\end{equation}
where
${\mathbb{X}}_t^{0,\dagger}=\lim_{N\rightarrow\infty}{X}_t^{0,\dagger},\
{\mathbb{X}}_t^{1,\dagger}=\lim_{N\rightarrow\infty}{X}_t^{1,\dagger},\  \overline{\alpha}_t^{\dag}=\lim_{N\rightarrow\infty}X_t^{(N),\dagger}$
are the associated limiting quantities, and $\mathbb{E}_t[\cdot]$ denotes the conditional expectation on the common information $\mathbb{F}^0$.
Similarly, taking the limit on \eqref{N-major-state}, we get the {asymptotic limit of the major agent}:
\begin{equation}\label{fu-4}
\left\{
\begin{aligned}
d\mathbb{X}_t^{0,\dagger}=&b^0(t,\mathbb{X}_t^{0,\dagger},u_t^0; {\overline{\alpha}}_t^{\dagger},\overline{u}_t^{*})dt
+\sigma^0(t,\mathbb{X}_t^{0,\dagger},u_t^0; {\overline{\alpha}}_t^{\dagger},\overline{u}_t^{*})dW_t^0,\quad \mathbb{X}_0^{0,\dagger}=x_0^0,\\
d{\overline{\alpha}}_t^{\dagger}=&
\mathbb{E}_t\Big[b(t,{\mathbb{X}}_t^{1,\dagger},u_t^{1,*};
\mathbb{X}_t^{0,\dagger},u_t^0;
{\overline{\alpha}}_t^{\dagger},
\overline{u}_t^{*})\Big]dt,\quad  \overline{\alpha}_0^{\dagger}=x_0,
\end{aligned}
\right.
\end{equation}
where ${\mathbb{X}}^{1,\dagger}$ denotes {the state of the representative minor agent, say, $\mathcal{A}_1$}
\begin{equation}\label{fu-5}
d{\mathbb{X}}_t^{1,\dagger}=b(t,{\mathbb{X}}_t^{1,\dagger},u_t^{1,*};
\mathbb{X}_t^{0,\dagger},u_t^0;{{\overline{\alpha}}_t^{\dagger}},
\overline{u}_t^{*})dt+\sigma(t,{\mathbb{X}}_t^{1,\dagger},u_t^{1,*};
\mathbb{X}_t^{0,\dagger},u_t^0;{{\overline{\alpha}}_t^{\dagger}},
\overline{u}_t^{*})dW_t^1,\quad
{\mathbb{X}}_0^{1,\dagger}=x_0.
\end{equation}
Comparing \eqref{fu-4} and \eqref{fu-5}, we have ${\overline{\alpha}}_t^{\dagger}=\mathbb{E}_t[{\mathbb{X}}^{1,\dagger}_t]$ by noting $\overline{u}_t^*=\mathbb{E}_t[u_t^{1,*}]$ and the solution uniqueness of the second equation in \eqref{fu-4}.
Then the limiting state $\mathbb{X}^{0,\dagger}$ of the major satisfies\begin{equation}\label{mfg-major-state}
\left\{
\begin{aligned}
d\mathbb{X}_t^{0,\dagger}=&b^0\big(t,\mathbb{X}_t^{0,\dagger},u_t^0; \mathbb{E}_t[{\mathbb{X}}^{1,\dagger}_t],\mathbb{E}_t[{u}_t^{1,*}]\big)dt
+\sigma^0\big(t,\mathbb{X}_t^{0,\dagger},u_t^0; \mathbb{E}_t[{\mathbb{X}}^{1,\dagger}_t],\mathbb{E}_t[{u}_t^{1,*}]\big)dW_t^0,\\
d{\mathbb{X}}_t^{1,\dagger}=&b\big(t,{\mathbb{X}}_t^{1,\dagger},
u_t^{1,*};
\mathbb{X}_t^{0,\dagger},u_t^0;\mathbb{E}_t[{\mathbb{X}}^{1,\dagger}_t],
\mathbb{E}_t[{u}_t^{1,*}]\big)dt
+\sigma\big(t,{\mathbb{X}}_t^{1,\dagger},
u_t^{1,*};
\mathbb{X}_t^{0,\dagger},u_t^0;\mathbb{E}_t[{\mathbb{X}}^{1,\dagger}_t],
\mathbb{E}_t[{u}_t^{1,*}]\big)dW_t^1,
\end{aligned}
\right.
\end{equation}with initial conditions $\mathbb{X}_0^{0,\dagger}=x_0^0,\ {\mathbb{X}}_0^{1,\dagger}=x_0.$
Taking limit on BSDE \eqref{fu-3} and similar to \eqref{mfg-major-state},
\begin{equation}\label{mfg-major-utility}
\left\{
\begin{aligned}
-d\mathbb{Y}_t^{0,\dagger}
=&f^0\Big(t,{\Pi_t^{0,\dagger},u_t^0;
\mathbb{E}_t[\Pi_t^{1,\dagger}],\overline{u}_t^{*}}
\Big)dt
-\mathbb{Z}_t^{0,\dagger}dW_t^0,\\
-d{\mathbb{Y}}_t^{1,\dagger}
=&f(t,{\Pi_t^{1,\dagger},u_t^{1,*};
\Pi_t^{0,\dagger},u_t^0;
\mathbb{E}_t[{u}_t^{1,*}],
\mathbb{E}_t[\Pi_t^{1,\dagger}]})dt
-{{\mathbb{Z}}_t^{1,\dagger,0}dW_t^0
-{\mathbb{Z}}_t^{1,\dagger,1}dW_t^1},\\
\mathbb{Y}_T^{0,\dagger}=&\Phi^0\big(\mathbb{X}_T^{0,\dagger},
\mathbb{E}_T[{\mathbb{X}}^{1,\dagger}_T]\big)+\xi^0,\quad {\mathbb{Y}}_T^{1,\dagger}=\Phi\big({\mathbb{X}}_T^{1,\dagger};
\mathbb{X}_T^{0,\dagger},\mathbb{E}_T[{\mathbb{X}}^{1,\dagger}_T]\big)+\xi^1,
\end{aligned}
\right.
\end{equation}
where $\Pi_t^{k,\dagger}=\lim_{N\rightarrow\infty}\Theta_t^{k,\dagger}
=(\mathbb{X}_t^{k,\dagger},\mathbb{Y}_t^{k,\dagger},\mathbb{Z}_t^{k,\dagger})$, $k=0,1.$
\begin{remark}\label{re-remainder}The remainder terms $\sum_{j=1}^{N}Z_t^{0,\dagger,j}dW_t^j$ and $\sum_{j=1, j \neq i}^{N}Z_t^{i,\dagger,j}dW_t^j$ in BSDE \eqref{fu-3} vanish as $N \longrightarrow +\infty$. For sake of presentation, we may omit these remainder terms hereafter.
\end{remark}
Finally, for the given $\mathbb{F}$-adapted ${u}^{1,*}$, we construct the following auxiliary problem for $\mathcal{A}_0$ associated with \eqref{mfg-major-utility} and \eqref{mfg-major-state}:
\begin{equation}\label{major-mfg1}
\begin{aligned}
 \sup_{u^0\in\mathcal{U}_0^c}\Big\{\Gamma^0(\mathbb{Y}_0^{0,\dagger})
 +\mathbb{E}\big[\int_0^Tg^0(t,\Pi_t^{0,\dagger},u_t^0;
\mathbb{E}_t[\Pi_t^{1,\dagger}],\overline{u}_t^{*})dt\big]\Big\}.\\
\end{aligned}
\end{equation}An optimal control (if exists) of  \eqref{major-mfg1} should rely on the given control $\big({u}^{1,*},\bar{u}_t^*\big)$.

\subsection{A bilateral perturbation: a representative minor agent}
We turn to a representative minor $\mathcal{A}_1$. By Definition \ref{def1}, $\mathcal{A}_1$ confronts an optimization problem when assuming $\mathcal{A}_0$ implements the exact Nash $u^{0,*}\in\mathcal{U}_0^c$; and $\mathcal{A}_{j} (j\geq 2)$ implement $u^{j,*}\in\mathcal{U}^c$.
If $\mathcal{A}_1$ applies a perturbed $u^1\in\mathcal{U}^c$, the state and functional of $\mathcal{A}_0, \mathcal{A}_j$ ($j\geq 2$) become
\begin{equation}\label{081301}
\left\{
\begin{aligned}dX_t^{j,*}=&b(t,X_t^{j,*},u_t^{j,*};
X_t^{0,*},u_t^{0,*};
X_t^{(N),*},u_t^{(N),*})dt+\sigma(t,X_t^{j,*},u_t^{j,*};
X_t^{0,*},u_t^{0,*};
X_t^{(N),*},u_t^{(N),*})dW_t^j,\\
-dY_t^{j,*}=&f(t,{\Theta_t^{j,*},u_t^{j,*};
\Theta_t^{0,*},u_t^{0,*};\Theta_t^{(N),*},u_t^{(N),*})dt
-Z_t^{j,*,0}dW_t^0-Z_t^{j,*,j}dW_t^j},\\
dX_t^{0,*}=&b^0(t,X_t^{0,*},u_t^{0,*}; {X_t^{(N),*},u_t^{(N),*}})dt+\sigma^0(t,X_t^{0,*},u_t^{0,*}; {X_t^{(N),*},u_t^{(N),*}})dW_t^0,\\
-dY_t^{0,*}=&f^0(t,\Theta_t^{0,*},u_t^{0,*};
\Theta_t^{(N),*},u_t^{(N),*})dt
-Z_t^{0,*}dW_t^0,
\end{aligned}
\right.
\end{equation}with $X_0^{0,*}=x_0^0,\ Y_T^{0,*}=\Phi^0(X_T^{0,*},X_T^{(N),*})+\xi^0,\ X_0^{j,*}=x_0,\ Y_T^{j,*}=\Phi(X_T^{j,*};X_T^{0,*},X_T^{(N),*})+\xi^j.$ We abuse notations (as with the same limits) to denote
$\Upsilon_\cdot^{(N),*}:=\frac{1}{N-1}\sum_{j= 2}^N\Upsilon_\cdot^{j,*},\ \Upsilon=X,Y,Z,u,$
and $\{\Theta_t^{k,*}=(X_t^{k,*},Y_t^{k,*},Z_t^{k,*})\}_{k=0}^{N};\ \ \Theta_t^{(N),*}=(X_t^{(N),*},Y_t^{(N),*},
Z_t^{(N),*}).$ Noting {$(\Theta^{0,*},\Theta^{j,*},\Theta^{(N),*})$} are
independent of $\mathcal{A}_1$, so they are \emph{exogenous} for $\mathcal{A}_1$. The state of $\mathcal{A}_1$ under $u^1$ satisfies
\begin{equation}\label{081401}
\begin{aligned}\nonumber
dX_t^{1,\ddagger}=&b(t,X_t^{1,\ddagger},u_t^1;X_t^{0,*},u_t^{0,*};
X_t^{(N),*},u_t^{(N),*})dt+\sigma(t,X_t^{1,\ddagger},u_t^1;X_t^{0,*},u_t^{0,*};
X_t^{(N),*},u_t^{(N),*})dW_t^1.
\end{aligned}
\end{equation}The \emph{quasi-realized} state $X_\cdot^{1,\ddagger}$ of $\mathcal{A}_1$ is with the superscript ``$\ddagger$" to indicate its dependence on perturbed $u^1$; whereas the \emph{exact-realized} state-control average $(X_\cdot^{(N),\ast},u_\cdot^{(N),\ast})$ only depend on the exact strategies so are still with ``$\ast$". The value $Y^{1,\ddagger}_0$} affected by $u^1$ of $\mathcal{A}_1$ satisfies\begin{equation}\label{052401}
\left\{
\begin{aligned}
-dY_t^{1,\ddagger}=&f(t,\Theta_t^{1,\ddagger},u_t^1;
\Theta_t^{0,*},u_t^{0,*};
\Theta_t^{(N),*},u_t^{(N),*})dt
-Z_t^{1,\ddagger,0}dW_t^0-Z_t^{1,\ddagger,1}dW_t^1,\\
Y_T^{1,\ddagger}=&\Phi(X_T^{1,\ddagger};X_T^{0,*},X_T^{(N),*})+\xi^1,
\end{aligned}
\right.
\end{equation}with $\Theta_t^{1,\ddagger}=(X_t^{1,\ddagger},Y_t^{1,\ddagger},Z_t^{1,\ddagger})$.
Next, similar to our analysis in Subsection 3.1, we can obtain the following coupled mean-field FBSDE from \eqref{081301} (noting $\mathbb{E}_t[u_t^{j,*}]=\overline{u}_t^*$), $j\geq2$,
 \begin{equation}\label{081701}
\left\{
\begin{aligned}
d\mathbb{X}_t^{0,*}=&b^0\big(t,\mathbb{X}_t^{0,*},u_t^{0,*}; \mathbb{E}_t[\mathbb{X}^{j,*}_t],\mathbb{E}_t[u_t^{j,*}]\big)dt
+\sigma^0\big(t,\mathbb{X}_t^{0,*},u_t^{0,*}; \mathbb{E}_t[\mathbb{X}^{j,*}_t],\mathbb{E}_t[u_t^{j,*}]\big)dW_t^0,\\
d\mathbb{X}_t^{j,*}=&b\big(t,\mathbb{X}_t^{j,*},u_t^{j,*};
\mathbb{X}_t^{0,*},u_t^{0,*};
\mathbb{E}_t[\mathbb{X}^{j,*}_t],\mathbb{E}_t[u_t^{j,*}]\big)dt
+\sigma\big(t,\mathbb{X}_t^{j,*},u_t^{j,*};
\mathbb{X}_t^{0,*},u_t^{0,*};
\mathbb{E}_t[\mathbb{X}^{j,*}_t],\mathbb{E}_t[u_t^{j,*}]\big)dW_t^j,\\
-d\mathbb{Y}_t^{0,*}=&f^0\Big(t,{\Pi_t^{0,*},u_t^{0,*};\mathbb{E}_t[\Pi^{j,*}_t],
\mathbb{E}_t[{u}^{j,*}_t]}\Big)dt
-\mathbb{Z}_t^{0,*}dW_t^0,\\
-d\mathbb{Y}_t^{j,*}=&f\Big(t,{\Pi_t^{j,*},u_t^{j,*};
\Pi_t^{0,*},u_t^{0,*};
\mathbb{E}_t[\Pi^{j,*}_t],
\mathbb{E}_t[{u}^{j,*}_t]\Big)dt
-\mathbb{Z}_t^{j,*,0}dW_t^0-\mathbb{Z}_t^{j,*,j}dW_t^j},\\
\mathbb{X}_0^{0,*}=&x_0,\ \mathbb{X}_0^{j,*}=\xi_j,\ \mathbb{Y}_T^{0,*}=\Phi^0\Big(\mathbb{X}_T^{0,*},
\mathbb{E}_T[\mathbb{X}^{j,*}_T]\Big)+\xi^0,\ \mathbb{Y}_T^{j,*}
=\Phi\Big(\mathbb{X}_T^{j,*};\mathbb{X}_T^{0,*},
\mathbb{E}_T[\mathbb{X}^{j,*}_T]\Big)+\xi^j,\\
\end{aligned}
\right.
\end{equation}by setting $\Pi_\cdot^{i,*}=(\mathbb{X}^{i,*}_\cdot,\mathbb{Y}^{i,*}_\cdot,\mathbb{Z}^{i,*}_\cdot)
=\lim_{N\rightarrow\infty}(X^{i,*}_\cdot,Y^{i,*}_\cdot,Z^{i,*}_\cdot),\
i \neq 1.$ Then by \eqref{081401}, \eqref{052401}, when $\mathcal{A}_1$ applying $u^1$, the limiting $\mathbb{X}^{1,\ddagger}(:=\lim_{N\rightarrow\infty}X^{1,\ddagger})$ and $\mathbb{Y}^{1,\ddagger}_t(:=\lim_{N\rightarrow\infty}Y_t^{1,\ddagger})$  satisfy
 \begin{equation}\label{081408}
\left\{
\begin{aligned}
d\mathbb{X}_t^{1,\ddagger}=&b\Big(t,\mathbb{X}_t^{1,\ddagger},u_t^1;
\mathbb{X}_t^{0,*},u_t^{0,*};
\mathbb{E}_t[\mathbb{X}^{j,*}_t],
\mathbb{E}_t[u_t^{j,*}]\Big)dt
+\sigma\Big(t,\mathbb{X}_t^{1,\ddagger},u_t^1;
\mathbb{X}_t^{0,*},u_t^{0,*};
\mathbb{E}_t[\mathbb{X}^{j,*}_t],
\mathbb{E}_t[u_t^{j,*}]\Big)dW_t^1,\\
-d\mathbb{Y}^{1,\ddagger}_t=&f\Big(t,{\Pi^{1,\ddagger}_t,u_t^1;
\Pi_t^{0,*},u_t^{0,*};
\mathbb{E}_t[\Pi^{j,*}_t],\mathbb{E}_t[{u}^{j,*}_t]}\Big)dt
{-\mathbb{Z}^{1,\ddagger,0}_tdW_t^0-\mathbb{Z}^{1,\ddagger,1}_tdW_t^1},\\
\mathbb{X}_0^{i,\ddagger}=&x_0,\quad
\mathbb{Y}^{1,\ddagger}_T=\Phi\Big(\mathbb{X}^{1,\ddagger}_T;
\mathbb{X}_T^{0,*},
\mathbb{E}_T[\mathbb{X}^{j,*}_T]\Big)+\xi^1,
\end{aligned}
\right.
\end{equation}where $\Pi_t^{1,\ddagger}=(\mathbb{X}_t^{1,\ddagger},\mathbb{Y}_t^{1,\ddagger},\mathbb{Z}_t^{1,\ddagger})$.
Along with  \eqref{081408}, the objective of  $\mathcal{A}_1$ is to maximize
\begin{equation}\label{080801}
\Gamma(\mathbb{Y}_0^{1,\ddagger})+\mathbb{E}\Big[\int_0^Tg(t,\Pi_t^{1,\ddagger},u_t^1;
\Pi_t^{0,*},u_t^{0,*}; \mathbb{E}_t[\Pi^{j,*}_t], \mathbb{E}_t[{u}^{j,*}_t])dt\Big].
\end{equation}Note that an optimal control of $\mathcal{A}_1$, if exists, should depend on ${u}^{0,*}$ and $\bar{u}^*$.

\subsection{A hierarchical recomposition}
We now apply a hierarchical recomposition to construct the desired auxiliary problem, that assumes a mixed two-layer game in the RMM context. To this end, we integrate the perturbed \eqref{mfg-major-state}-\eqref{major-mfg1} in side of $\mathcal{A}_0$ indexed by $``\dag"$, and \eqref{081408}-\eqref{080801} in side of the representative minor $\mathcal{A}_1$ by $``\ddag"$ together. This yields an (extended) state $({X}^{0,\dagger},{X}^{1,\dagger};{X}^{{1,}\ddagger},{X}^{0,\ddagger},{X}^{j,\ddagger})$ labeled by $j\geq 2$:
\begin{equation}\label{mlf-state}
\left\{
\begin{aligned}
\quad d{X}_t^{0,\dagger}
=&b^0\big(t,{X}_t^{0,\dagger},
{{u_t^0}};
\mathbb{E}_t[{{X}}^{1,\dagger}_t],
\mathbb{E}_t[{u_t^{1}}]\big)dt
+\sigma^0\big(t,{X}_t^{0,\dagger},
{{u_t^0}};
 \mathbb{E}_t[{{X}}^{1,\dagger}_t],
\mathbb{E}_t[{u_t^{1}}]\big)dW^0_t,\\
d{{X}}_t^{1,\dagger}=&b\big(t,{{X}}_t^{1,\dagger},
{u_t^{1}};
{X}_t^{0,\dagger},{u_t^0};
\mathbb{E}_t[{{X}}^{1,\dagger}_t],
\mathbb{E}_t[{u_t^{1}}]\big)dt
+\sigma\big(t,{{X}}_t^{1,\dagger},
{u_t^{1}};
{X}_t^{0,\dagger},{u_t^0};
\mathbb{E}_t[{{X}}^{1,\dagger}_t],
\mathbb{E}_t[{u_t^{1}}]\big)dW_t^1,\\
\quad d{X}_t^{{1,}\ddagger}
=&b\big(t,{X}_t^{{1},\ddagger},
{u_t^1};
{{X}_t^{0,\ddagger},{u_t^{0}};
\mathbb{E}_t[{X}^{j,\ddagger}_t],
\mathbb{E}_t[{u_t^{j}}]}\big)dt+
\sigma\big(t,{X}_t^{{1},\ddagger},
{u_t^1};
{{X}_t^{0,\ddagger},{u_t^{0}};
\mathbb{E}_t[{X}^{j,\ddagger}_t],
\mathbb{E}_t[{u_t^{j}}]}\big)dW_t^1,\\
\quad d{X}_t^{0,\ddagger}=&b^0\big(t,{X}_t^{0,\ddagger},
{u_t^{0}}; \mathbb{E}_t[{X}^{j,\ddagger}_t],
\mathbb{E}_t[{u_t^{j}}]\big)dt
+b^0\big(t,{X}_t^{0,\ddagger},
{u_t^{0}}; \mathbb{E}_t[{X}^{j,\ddagger}_t],
\mathbb{E}_t[{u_t^{j}}]\big)dW_t^0,\\
\quad d{X}_t^{j,\ddagger}=&b\big(t,{X}_t^{j,\ddagger},
{u_t^{j}};
{X}_t^{0,\ddagger},{u_t^{0}};
\mathbb{E}_t[{X}^{j,\ddagger}_t],
\mathbb{E}_t[{u_t^{j}}]\big)dt
+
\sigma\big(t,{X}_t^{j,\ddagger},
{u_t^{j}};
{X}_t^{0,\ddagger},{u_t^{0}};
\mathbb{E}_t[{X}^{j,\ddagger}_t],
\mathbb{E}_t[{u_t^{j}}]\big)dW_t^j,\\
{X}_0^{0,\dagger}=&{X}_0^{0,\ddagger}=x^0_0,\quad {{X}}_0^{1,\dagger}={X}_0^{1,\ddagger}={X}_0^{j,\ddagger}=x_0.
\end{aligned}
\right.
\end{equation}Indeed, the first and second state $({X}^{0,\dagger}, {{X}}^{1,\dagger})$ come from \eqref{mfg-major-state}, the third ${X}^{1,\ddagger}$ from \eqref{081408}, and the last two $({X}^{0,\ddagger},{X}^{j,\ddagger})$ from  \eqref{081701} by replacing exact centralized controls $(u^{0,*}, u^{1,*}, u^{j,*})$ with decentralized ones $(u^{0},u^1,u^{j}).$ With the extended state \eqref{mlf-state}, there arise three agents, respectively: a follower $\mathcal{A}_{0}$ using decentralized control $u^{0}$, a follower $\mathcal{A}_{1}$ using $u^1$, and a leader $\mathcal{A}_{j}$ using $u^{j},$ for some $j\geq 2.$ Recall that $u^0$ is $\mathbb{F}^0$-adapted and $u^1$ (resp. $u^j$) is $\mathbb{F}^{1,0}$- (resp. $\mathbb{F}^{j,0}$)-adapted. Therefore, the follower $\mathcal{A}_{0}$ can \emph{directly} affect all 5-tuple components, the  \emph{follower} $\mathcal{A}_{1}$ \emph{directly} affects the first three $({X}^{0,\dagger}, {{X}}^{1,\dagger},{X}^{1,\ddagger}),$ while
the  \emph{leader} $\mathcal{A}_{j}$ \emph{directly} affect the last three $({X}^{1,\ddagger}, {X}^{0,\ddagger}, {X}^{j,\ddagger}).$ On the other hand, all 5-tuple components are coupled through their dependence on the control triple $(u^{0}, u^1, u^{j}).$ In this sense, $\mathcal{A}_{1}$ and $\mathcal{A}_{j}$ also influence the 5-tuple (\emph{indirectly}); all components are endogenous and not redundant.

We can now formulate limiting recursive functionals. Specifically, the follower $\mathcal{A}_{0}$ aims to
\begin{equation}\label{feb2901}
 {\sup\limits_{u^0\in\mathcal{U}^d_0}}\Big\{\Gamma^0({Y}_0^{0,\dagger})
 +\mathbb{E}\big[\int_0^Tg^0(t,\Theta_t^{0,\dagger},u_t^0;
 \mathbb{E}_t[{\Theta}_t^{1,\dagger}], \mathbb{E}_t[u_t^1])dt\big]\Big\},
\end{equation}
where $({Y}^{0,\dagger},{Z}^{0,\dagger};{Y}^{1,\dagger},{Z}^{1,\dagger})$ is the solution of
the coupled mean-field BSDE
\begin{equation}\label{mlf-major-utility}
\left\{
\begin{aligned}
-d{Y}_t^{0,\dagger}=&f^0\Big(t,{\Theta_t^{0,\dagger},u_t^0;
\mathbb{E}_t[{\Theta}_t^{1,\dagger}],\mathbb{E}_t[{u}_t^{1}]
}\Big)dt
-{Z}_t^{0,\dagger}dW_t^0,\\
-d{{Y}}_t^{1,\dagger}=&f(t,
{\Theta_t^{1,\dagger},u_t^{1};
\Theta_t^{0,\dagger},u_t^0;
\mathbb{E}_t[\Theta_t^{1,\dagger}],
\mathbb{E}_t[{u}_t^{1}]
)dt
-{{Z}}_t^{1,\dagger,0}dW_t^0-{{Z}}_t^{1,\dagger,1}dW_t^1},\\
{Y}_T^{0,\dagger}=&\Phi^0\big({X}_T^{0,\dagger},
\mathbb{E}_T[{{X}}^{1,\dagger}_T]\big)+\xi^0,\quad {{Y}}_T^{1,\dagger}=\Phi\big({{X}}_T^{1,\dagger};
{X}_T^{0,\dagger},\mathbb{E}_T[{{X}}^{1,\dagger}_T]\big)+\xi^1.
\end{aligned}
\right.
\end{equation}
The aim of the follower $\mathcal{A}_{1}$ is
\begin{equation}\label{feb2902}
{\sup_{u^1\in\mathcal{U}_1^d}}
\Big\{\Gamma({Y}_0^{1,\ddagger})
+\mathbb{E}\Big[\int_0^Tg(t,\Theta_t^{1,\ddagger},u_t^1;
\Theta_t^{0,\ddagger},u_t^{0};
\mathbb{E}_t[\Theta^{j,\ddagger}_t],\mathbb{E}_t[{u}^{j}_t])dt\Big]\Big\},
\end{equation}
where $({Y}^{1,\ddagger},{Z}^{1,\ddagger})$ is the solution of the BSDE
 \begin{equation}\label{081411-con}
\left\{
\begin{aligned}
-d{Y}^{1,\ddagger}_t=&f\Big(t,
\Theta^{1,\ddagger}_t,u_t^1;
\Theta_t^{0,\ddagger},u_t^{0};
\mathbb{E}_t[\Theta^{j,\ddagger}_t],\mathbb{E}_t[{u}^{j}_t]\Big)dt
-{Z}^{1,\ddagger,0}_tdW_t^0-{Z}^{1,\ddagger,1}_tdW_t^1,\\
{Y}^{1,\ddagger}_T=&\Phi\Big({X}^{1,\ddagger}_T;
{X}_T^{0,\ddagger},
\mathbb{E}_T[{X}^{j,\ddagger}_T]\Big)+\xi^1,
\end{aligned}
\right.
\end{equation}
with the exogenous processes
$({Y}^{0,\ddagger},{Z}^{0,\ddagger};{Y}^{j,\ddagger},{Z}^{j,\ddagger})$ given by
\begin{equation}\label{081702-con}
\left\{
\begin{aligned}
-d{Y}_t^{0,\ddagger}=&f^0\Big(t,\Theta_t^{0,\ddagger},u_t^{0};
\mathbb{E}_t[\Theta^{j,\ddagger}_t],\mathbb{E}_t[{u}^{j}_t]\Big)dt
-{Z}_t^{0,\ddagger}dW_t^0,\\
-d{Y}_t^{j,\ddagger}=&f\Big(t,\Theta_t^{j,\ddagger},u_t^{j};
\Theta_t^{0,\ddagger},u_t^{0};
\mathbb{E}_t[\Theta^{j,\ddagger}_t],\mathbb{E}_t[{u}^{j}_t]\Big)dt
-{Z}_t^{j,\ddagger,0}dW_t^0-{Z}_t^{j,\ddagger,j}dW_t^j,\\
{Y}_T^{0,\ddagger}=&\Phi^0\Big({X}_T^{0,\ddagger},
\mathbb{E}_T[{X}^{j,\ddagger}_T]\Big)+\xi^0,\quad
{Y}_T^{j,\ddagger}
=\Phi\Big({X}_T^{j,\ddagger};{X}_T^{0,\ddagger},
\mathbb{E}_T[{X}^{j,\ddagger}_T]\Big)+\xi^j.\\
\end{aligned}
\right.
\end{equation}
And the leader $\mathcal{A}_{j}$ aims to minimize the following \emph{quadratic deviation} functional
\begin{equation}\label{320}
{\inf_{u^{j}\in\mathcal{U}^d_j}||u^{j}-u^{1}||_{L^{2}}^{W^0}:=
\mathbb{E}\int_0^{T}|\mathbb{E}_{t}[u^{j}_{t}]-\mathbb{E}_{t}[u^{1}_{t}]|^{2}dt.}
\end{equation}
Last, we identify mixed leader-follower-Nash interactions among $\mathcal{A}_0, \mathcal{A}_1, \mathcal{A}_j$ as follows:
\begin{itemize}
\item As the unique leader, at the beginning $t=0,$ $\mathcal{A}_{j}$ (for some $j\geq 2$) announces an (open-loop) control $u^{j}$ on the whole horizon $[0, T]$ to the followers $\mathcal{A}_{0}, \mathcal{A}_{1}$ both.
\item Given the announced $u^{j}|_{[0, T]},$ $\mathcal{A}_{0}, \mathcal{A}_{1}$ will find a Nash equilibrium $u^{0}, u^{1}$ simultaneously. In particular, $\mathcal{A}_{0}$ and $\mathcal{A}_{1}$ optimize the functionals \eqref{feb2901} and \eqref{feb2902}, respectively.

    \item Anticipate the best response $u^{0}$ and $u^{1}$ parameterized by the priori $u^{j}$, the leader intakes such leader-follower interaction to minimize the quadratic deviation cost. Such minimum can be reached via a fixed point argument or consistency condition.
\end{itemize}Given the above interactions, a mixed triple-agent game should be solved by the following steps.\\
(1) First, to formulate and solve an auxiliary control problem for  $\mathcal{A}_0$ in terms of control $u^{0}$ and state $({X}^{0,\dagger}, {{X}}^{1,\dagger},{Y}^{0,\dagger},{Z}^{0,\dagger},{Y}^{1,\dagger},{Z}^{1,\dagger})$ by fixing the generic admissible control $u^1$ of $\mathcal{A}_{1}$. We may denote the optimal one as $\tilde{u}^{0}=\tilde{u}^{0}[u^1]$ to show its dependence on $u^1$. \\
(2) Second, for given ${u}^{0}$ and pre-announced $u^{j}$, to formulate and solve an auxiliary control problem for  $\mathcal{A}_1$ in terms of control $u^{1}$ and state $({X}^{1,\ddagger}, {X}^{0,\ddagger}, {X}^{j,\ddagger},{Y}^{1,\ddagger},{Z}^{1,\ddagger},{Y}^{0,\ddagger},{Z}^{0,\ddagger},$ ${Y}^{j,\ddagger},{Z}^{j,\ddagger}).$
We may denote the optimal one as $\tilde{u}^{1}=\tilde{u}^{1}[u^{j},{u}^{0}]$ depending on $u^j$ and $u^0.$ \\
(3) Last, get the Nash equilibrium of $\mathcal{A}_0$ and $\mathcal{A}_1$, say $(\hat{u}^{0}[u^{j}],\hat{u}^{1}[u^{j}])$ depending on  pre-announced $u^j$; and solve the optimization of the leader by matching
$\mathbb{E}_t[u^{j}_t]
=\mathbb{E}_t\big[\hat{u}_t^{1}[u^{j}_t]\big].$

\begin{remark} (i) When the weak-coupling of RMM does not include the control-average of all minors, the term $\mathbb{E}_t[{X}^{j}_t]$ in the third and fifth equation of \eqref{mlf-state} can be replaced by some off-line process. The last equation in \eqref{mlf-state} thus becomes redundant. In fact, the inclusion of the control-average of all minors necessitates the introduction of additional dynamics to account for the averages on the physical, recursive and intensity state.

(ii) Furthermore, if we exclude the consideration of recursive functionals (namely, $\Gamma^0=\Gamma\equiv 0$), the mixed triple-agent game will simplify to a two-agent nonzero-sum game, as in \cite{CZ2016}.

(iii) The introduction of the additional $\mathcal{A}_j$ facilitates a clearer specification of the exogenous processes, as discussed in Subsection 3.2 from the standpoint of $\mathcal{A}_1$. Otherwise, two exogenous processes s would need to be introduced to replace $(\mathbb{E}_t[\mathbb{X}_t^{j,*}],\mathbb{E}_t[u_t^{j,*}])$ simultaneously, accompanied by an associated fixed point analysis. The overall analysis will thus become complicated, given that the consistency matching already involves another fixed-point analysis. In contrast, our formulation of the mixed triple-agent game, particularly the introduction of the virtual $\mathcal{A}_{j}$, may clearly elucidate the intricate exogenous-endogenous relationships.\end{remark}

\subsection{A unified structural scheme}

Our bilateral perturbations and mixed hierarchical recomposition in Sections 3.2 and 3.3 indeed suggest a unified scheme to analyze general LP systems with more complex couplings.

We start with a generic LP system characterized by a key 5-tuple $(\mathcal{I}, \mathcal{A}, J, \mathcal{D}, \mathcal{C}_{0})$ where (i) $\mathcal{I}$ is the index set of all agents, as enumerated by $\mathcal{A}:=\{\mathcal{A}_{i}\}_{i \in \mathcal{I}}$; (ii) $J:=\{J_{i}\}_{i \in \mathcal{I}}$ the set of weak-coupled functionals to be optimized by $\mathcal{A}_{i}, i \in \mathcal{I},$ respectively; (iii) $\mathcal{D}:=\{u_{i}\}_{i \in \mathcal{I}}$ denotes associated decision profile of $\mathcal{A}$; (iv)
$\mathcal{C}_{0}:=\{\mathcal{I}_{k}\}_{k \in \Theta}$ is the \emph{maximal} coalition structure on $\mathcal{I}$, which is structurally determined by $J, \mathcal{D}$ and will be elaborated later. Typically, $\mathcal{I}=\mathcal{I}_1:=\{1, \cdots, N\}$, $\mathcal{I}_2:=\{1, \cdots, k\} \times \mathcal{I}_1$, or $\mathcal{I}_3:=\{1, \cdots, k\} \cup \mathcal{I}_1$ respectively, represent the classical LP systems with homogenous agents, $k$ classes of heterogenous agents, or $k$ major agents with homogenous minors. For second order countable or continuum infinity, we can let $\mathcal{I}_4:=\mathcal{I}_1 \times \mathcal{I}_1$ or $\mathcal{I}_5=[0, 1].$ Given $(\mathcal{I}, J, \mathcal{D}, \mathcal{C}_{0})$, a structural scheme includes the following three steps.

Step 1: \textbf{Exchangeable decomposition}. Essentially, the MFG is an effective dimension-reduction analysis relying on exchangeabilities across agents, shares a similar spirit with the notable \emph{symmetric game} in deterministic context. Exchangeabilities of a controlled $(\mathcal{I}, \mathcal{A}, J, \mathcal{D}, \mathcal{C}_{0})$ can be characterized asymptotically, as $|\mathcal{I}| \longrightarrow +\infty$, by the so-called \emph{coalition structure}. This refers to a partition of the index set $\mathcal{I}$ (namely, $\sum_{k \in \Theta} \mathcal{I}_{k}=\mathcal{I}$ where $\sum$ the disjoint sets' union); and for $\forall k \in \Theta$, all agents $\{\mathcal{A}_{j}| j \in \mathcal{I}_{k}\}$ form an exchangeable sub-class, in the sense that\begin{equation}
\begin{aligned}\nonumber
&(\mathcal{I}, J, \mathcal{D}, \mathcal{C}_{0})=(\sum_{k \in \Theta}\mathcal{I}_{k}, \ \ \sum_{j \in \mathcal{I}_{k}}\sum_{k \in \Theta}J_{j,k}, \ \ \sum_{j \in \mathcal{I}_{k}}\sum_{k \in \Theta}u_{j,k},  \ \  \mathcal{C}_{0})\\&=(\mathcal{I}_{k} \cup \sum_{k'\neq k \in \Theta}\mathcal{I}_{k'}, \ \ \{J_{j,k}\}_{j \in \mathcal{I}_{k}} \cup \sum_{j \in \mathcal{I}_{k'}}\sum_{k' \neq k \in \Theta}J_{j,k'}, \ \ \{u_{j,k}\}_{j \in \mathcal{I}_{k}} \cup \sum_{j \in \mathcal{I}_{k'}}\sum_{k'\neq k \in \Theta}u_{j,k'}, \ \ \mathcal{C}_{0})\\& \Longleftrightarrow_{\mathcal{E}} (\widetilde{\mathcal{I}_{k}} \cup \sum_{k'\neq k \in \Theta}\mathcal{I}_{k'}, \ \ \ \widetilde{\{J_{j,k}\}}_{j \in \mathcal{I}_{k}} \cup \sum_{j \in \mathcal{I}_{k'}}\sum_{k' \neq k \in \Theta}J_{j,k'}, \ \ \ \widetilde{\{u_{j,k}\}}_{j \in \mathcal{I}_{k}} \cup \sum_{j \in \mathcal{I}_{k'}}\sum_{k'\neq k \in \Theta}u_{j,k'}, \ \ \ \mathcal{C}_{0}) \\& := (\widetilde{\mathcal{I}}, \widetilde{J}, \widetilde{\mathcal{D}}, \mathcal{C}_{0})
\end{aligned}
\end{equation}where $\left(\widetilde{\mathcal{I}_{k}}, \widetilde{\{J_{j,k}\}}_{j \in \mathcal{I}_{k}}, \widetilde{\{u_{j,k}\}}_{j \in \mathcal{I}_{k}}\right)$ is a simultaneous permutation on the sub-index set $\mathcal{I}_{k}$ and associated intersections on $J$ and $\mathcal{D}$; $``\Longleftrightarrow_{\mathcal{E}}"$ denotes the equivalence relation in terms of game equilibrium. That is, the set of Nash equilibria of $(\mathcal{I}, J, \mathcal{D}, \mathcal{C}_{0})$ is invariant to that of $(\widetilde{\mathcal{I}}, \widetilde{J}, \widetilde{\mathcal{D}}, \mathcal{C}_{0})$ under finite permutation on $\mathcal{I}_{k}$. Roughly, this means that the simultaneous optimizations faced by $\{\mathcal{A}_{j}\}_{j \in \mathcal{I}_{k}}$, are endowed with the identical probabilistic structures. Moreover, under large symmetric assumption as in \cite{CZ2018},
the equivalence $\mathcal{E}$ should be element-wise and thus transitive. Therefore, there exists a \emph{maximal} coalition, denoted as $\mathcal{C}_{0}$, provided the set of coalition structures is non-empty. By ``maximal", it is the coarsest partition than other coalitions, that implies the largest exchangeable decomposition and hence the largest dimension reduction can be achieved. In fact, this maximal coalition can be constructed through the \emph{saturated sets} using equivalent relation, as the union of all $\mathcal{E}$-equivalent classes (\cite{HL2019}).

A trivial coalition is $\{\{e\} | e \in \mathcal{I}\},$ the set of all singleton sets generated by each element. In this case, $|\Theta|=|\mathcal{I}|$ without dimension reduction. The MFG is applicable for non-trivial coalitions in that $|\Theta|=O(1)$ as $|\mathcal{I}|\longrightarrow +\infty$; or reversely, for at least one $k \in \Theta,$ $|\mathcal{I}_{k}|=O(|\mathcal{I}|).$ For example, $|\mathcal{I}_{k}|=O(|\mathcal{I}|^{\alpha})$ for some $\alpha >0$. For LP systems with $\mathcal{I}_{1,2,3},$ $\mathcal{C}_0=\mathcal{I}_{1}, (\{1\}\times \mathcal{I}_1, \cdots, \{k\}\times \mathcal{I}_1), (\{1\}, \cdots, \{k\}, \mathcal{I}_1),$ with associated $|\Theta|=1, k, k+1$, respectively. For these cases, heuristic arguments are still feasible to construct auxiliary control. However, they boil down when coupling structures assume more complex forms, such as non-standard bridge configuration; and in this case, the identification of $\mathcal{C}_0$ with associated decomposition become necessary to alternatively yield a unified and systematic analysis.

Step 2: \textbf{Representative multilateral perturbation}. Given the maximal coalition $\mathcal{C}_{0}$, one can select a representative agent, denoted as $\mathcal{A}^{\text{rep}}_{k}, k \in \Theta$, from each exchangeable sub-class $\{\mathcal{A}_{j}| j \in \mathcal{I}_{k}\}.$ This selection yields a representative collection $\mathcal{R}:=\{\mathcal{A}^{\text{rep}}_{k}\}_{k \in \Theta}$ that is dimension-reduced by noting $|\Theta|=O(1)$. In fact, we may abuse notations to denote $\mathcal{R}={\mathcal{A}}/{\Theta}$ as the quotient space with elements in the equivalent classes by $\mathcal{E}$. Then, a multilateral perturbation can be introduced in side of each $\mathcal{A}^{\text{rep}}_{k}$ separately, by assuming all other representatives still keep their equilibrium strategies. Depending on the weak-coupling mechanism structured by $\{{J_{j}}\}_{j \in \mathcal{I}}$, each perturbation $u^{\text{rep}}_{k}$ will be transmitted throughout the whole LP system $\{\mathcal{A}_{j}\}_{j \in \mathcal{I}}=\sum_{k \in \Theta}\{{\mathcal{A}_{j}}\}_{j \in \mathcal{I}_{k}}$ across all exchangeable sub-classes. A typical transmission, in an open-loop and dynamic setting, can be sketched by the following channel via the weak-coupling of state-average:$$u^{\text{rep}}_{k}\Longrightarrow x^{\text{rep}}_{k}(u^{\text{rep}}_{k}) \Longrightarrow \left(x^{(N)}_{k}(x^{\text{rep}}_{k}(u^{\text{rep}}_{k})), \ \ x^{(N)}_{k' \neq k}(x^{\text{rep}}_{k}(u^{\text{rep}}_{k}))\right) \Longrightarrow J_{j' \neq j}(x^{(N)}_{k}, x^{(N)}_{k' \neq k}) \Longrightarrow \cdots $$Along with such transmission, influence of a representative $\mathcal{A}^{\text{rep}}_{k}$ on controlled LP systems can be completely quantified, that is essentially equivalent to the Fr\'{e}chet differential of $J_{k}^{\text{rep}}$ on $u^{\text{rep}}_{k}$.

Step 3: \textbf{Hierarchical recomposition.}
The multilateral perturbation involves approximations to completely quantify all LP weak-couplings asymptotically, thereby leading to a variety of limiting quantities. Due to the exchangeability, these quantities typically assume relevant conditional expectations on the tail sigma-algebra, as per {DeFinetti theorem}. However, unlike those in classical McKean-Vlasov control problems, these expectations exhibit distinct modes in realized degrees of controllability, contingent upon their hierarchical positions within the entire weak-coupled structure. For instance, in the RMM, three modes emerge: exact-realized (exact limit), quasi-realized (semi-exact), and null-realized, as indicated by \eqref{mfg-major-state} and \eqref{081301}.
In the current RMM or simpler coupled setups, these modes can be ordered using an {exogenous-endogenous} relation, where the ``exogenous" exact-realized mode dominates the ``endogenous" null-realized one. However, in more complex forms of LP couplings, these modes cannot be fully encapsulated by a binary ordering alone, instead forming a more intricate directed graph network, which is challenging to be studied by heuristic arguments alone.

Our resolution is a hierarchical recomposition, similar to the well-studied structural function through the so-called \emph{path sets} in reliability theory (\cite{AB1984}).
In fact, reliability of any system is equivalent to that of a serial (sequential) arrangement of parallel subsystems. Likewise, for a generic LP system, all involved modes by Step 2 can be stratified into a sequential layers with leader-follower-type hierarchies; and each layer consists of parallel Nash-type nodes with simultaneous decisions. This stratification, similar to the construction of structural function, is indeed applicable to any LP systems, such as those with \emph{non-classical intermediate} coupling, akin to the \emph{bridge systems} in reliability analysis. Such stratification enables an unified and more clear-cut construction of the desired auxiliary problem by recomposing all modes across hierarchical layers with a combination of corresponding variants of fixed-point matching. For example, the aforementioned three modes in RMM yield two layers, thereby enabling the construction of the auxiliary problem via a triple-agent leader-follower-Nash game. By coincidence, the multilateral perturbation assumes a similar role to the synthesis of all \emph{minimal path sets }in reliability analysis, as both aim to quantify all transmission channels within a given system.

\textbf{Attainability}. For classical LP systems, leader-type agents at higher layers in the hierarchical recomposition often engage in a fixed-point analysis, as demonstrated by \eqref{320} for the RMM when $\mathcal{A}_{j}$ aligns the announced $u^{j}$ with the resultant $u^{1}(u^{j})$ on their conditional expectations $\mathbb{E}_{t}u^{j}, \mathbb{E}_{t}u^{1}$. In contrast to follower-type agents at lower layers, who confront more regular control problems (see \eqref{feb2902}), the fixed-point analysis for leaders can indeed degenerate, as the minimal deviations become trivially attainable at zero, provided $\mathbb{E}_{t}u^{j}=\mathbb{E}_{t}u^{1}$. Consequently, the choice of norm $\phi(\cdot)$ on deviation $\mathbb{E}_{t}u^{j}-\mathbb{E}_{t}u^{1}(u_{j})$ becomes indifferent. Specifically, the introduction of a quadratic deviation on $L^{2}$-norm: {$\phi(\cdot)=||\cdot||_{L^{2}}^{W^{0}}$} as in \eqref{320}, is merely formal. In fact, the top layer, characterized by a high degree of realized controllability, presents a tradeoff where its remaining control capability is more prone to degeneration. However, in LP systems with more complex couplings, there may emerge nonclassical hierarchical layers and the relevant analysis, particularly of top layers, may no longer be degenerate. This is particularly true for LP systems with \emph{nested} or \emph{asymmetric} information (\cite{AMS1995,K2019}), 
or those with heterogenous robustness beliefs, as well as those with bridging-intermediate-type couplings. In such instances, the optimal deviation norms cannot be trivially attainable at zero, necessitating the replacement of fixed-point analysis with some non-trivial optimization problems.

\textbf{Summary}. Step 1-3 constitute an unified structural scheme to analyze more general LP systems, particularly those with non-classical coupling structures. In fact, for classical LP systems with only minor (homogenous or heterogenous) agents, auxiliary problems can be constructed by straightforward heuristic arguments. The main reason is that their coalition structure are relatively simple thus no need to invoke Step 2-3. A more complex but classical LP system is the one involving single major agent, for which auxiliary problem can still be constructed heuristically but shall invoke an exogenous-endogenous relation to tackle additional couplings by the major agent. Our RMM has fairly general weak-couplings, especially those of the recursive state pairs, that motivate us to introduce the structural scheme, especially Step 2-3, to fully quantify complexities of all resultant perturbation transmissions and construct the auxiliary problem hierarchically. For non-classical LP systems such as those with nested or asymmetric information, or those featuring bridge-type couplings, heuristic arguments are no longer feasible for constructing auxiliary problems. For example, even for LP systems with $k$ majors and $K$ classes of heterogenous minors, heuristic construction necessitates an exogenous-endogenous structure represented by a directed graph with $k \times K$ edges. By contrast, the structural scheme, conversely, simplifies this to a sequential flow with only $k+K$ binary relations. 

\section{Approximated Nash equilibrium of the RMM problem}

\subsection{Auxiliary control problem for the follower  $\mathcal{A}_0$}
Now we fix the follower $\mathcal{A}_1$'s control $u^1$ and solve the auxiliary problem of $\mathcal{A}_0$.
 The auxiliary problem for $\mathcal{A}_0$ is formulated by \eqref{feb2901} 
associated with mean-field BSDE \eqref{mlf-major-utility} and McKean-Vlasov SDE \eqref{mlf-state} (its first two components).
Related Hamiltonian is defined as$$H_0
=\langle p^0,b^0\rangle
+\langle p,b\rangle
+\langle q^{00},\sigma^0\rangle
+\langle q^{11},\sigma\rangle
- l^0\cdot f^0
- l\cdot f
+g^0,$$where $(b^0,\sigma^0)=(b^0,\sigma^0)(t,x_0,u_0;\bar{x},\bar{u}),$
$(b,\sigma)=(b,\sigma)(t,x,u;x_0,u_0;\bar{x},\bar{u})$,
$(f^0,g^0)=(f^0,g^0)(t,x_0,$ $y_0,z_0,u_0;\bar{x},\bar{y},\bar{z},$ $\bar{u})$,
$f=f(t,x,y,z,u;x_0,y_0,z_0,u_0;\bar{x},\bar{y},\bar{z},\bar{u})$.
Using stochastic maximum principle (\cite{AD2011}), we get the following Hamiltonian system
\begin{equation}\label{H-system-major}
\left\{
\begin{aligned}
d{X}_t^{0,*}=&\partial_{p_0}H_0(t)dt
+\partial_{q_{00}}H_0(t)dW_t^0,\
d{{X}}_t^{1,*}=\partial_{p}H_0(t)dt
+\partial_{q_{11}}H_0(t)dW_t^1, \\
-d{Y}_t^{0,*}=&-\partial_{l^0}H_0(t)dt
-{Z}_t^{0,*}dW_t^0,\
-d{{Y}}_t^{1,*}=-\partial_{l}H_0(t)dt
-{{Z}}_t^{1,*,0}dW_t^0-{{Z}}_t^{1,*,1}dW_t^1, \\
dL_t^0=&-\partial_{y_0}H_0(t)dt-\partial_{z_0}H_0(t)dW^0_t,\\
dL_t=&-\big[\partial_{y}H_0(t)
+\mathbb{E}_t[\partial_{\bar{y}}H_0(t)]\big]dt
-\big[\partial_{z}H_0(t)+\mathbb{E}_t[\partial_{\bar{z}}H_0(t)]\big]d(W^0,W^1)_t,
\\
-dP_t^0=&\partial_{x_0}H_0(t)dt-Q_t^{00}dW^0_t-Q_t^{01}dW_t^1, \\
-dP_t=&\big[\partial_{x}H_0(t)
+\mathbb{E}_t[\partial_{\bar{x}}H_0(t)]\big]dt-Q_t^{10}dW^0_t-Q_t^{11}dW_t^1,\\
\end{aligned}
\right.
\end{equation}
where ${X}_0^{0,*}=x_0^0$,  $X_0^{1,*}=x_0,$ ${Y}_T^{0,*}=\Phi^0\big({X}_T^{0,*},
\mathbb{E}_T[{{X}}^{1,*}_T]\big)+\xi^0$, $ Y_T^{1,*}=\Phi\big({{X}}_T^{1,*};
{X}_T^{0,*},\mathbb{E}_T[{{X}}^{1,*}_T]\big)+\xi^1,$ $L_0^0=-\Gamma^0_y(Y_0^{0,*}),$ $  L_0=0$, $P_T^0=-L_T^0\cdot\partial_{x_0}\Phi^0
-L_T\cdot\partial_{x_0}\Phi,$ $P_T=-L_T\cdot\partial_x\Phi
-\mathbb{E}_T[L_T^0\cdot\partial_{\bar{x}}\Phi^0+L_T\cdot\partial_{\bar{x}}\Phi]$; $(L^0,L;P^0,P,Q^{00},Q^{01},Q^{10},Q^{11})$ are the associated adjoint processes and
\begin{equation}\label{070601}
H_0(t):=H_0(t,\Theta_t^{0,*},u_t^{0,*};
\Theta_t^{1,*},u_t^1;
\mathbb{E}_t[\Theta_t^{1,*}],\mathbb{E}_t[u_t^1];
L_t^0,L_t;P_t^0,P_t,Q_t^{00},Q_t^{11}),
\end{equation}
assuming that $u^{0,*}$ is an optimal control and {$(\Theta^{0,*},\Theta^{1,*})=(X^{0,*},Y^{0,*},Z^{0,*};X^{1,*},Y^{1,*},Z^{1,*})$} the associated  optimal state trajectory.
Then we have the following result.
\begin{proposition}\label{pro1}
Let  (A1) be in force. Moreover, we assume that \\
(1) There exists a unique maximizer of the Hamiltonian $H_0$ as a function of $u_0$
(denoted by $\widehat{u}^0$);\\
(2) The function $H_0$ is convex in $(x_0,x,y_0,y,z_0,z,u_0)$.\\
If $(\Theta^{0,*},\Theta^{1,*};L^0,L;P^0,Q^{00},Q^{01};P,Q^{10},Q^{11})$ solves system \eqref{H-system-major}, the optimal control of $\mathcal{A}_0$ is
$$\widetilde{u}_t^0=\mathbb{E}_t\big[
\widehat{u}^0\big(t,\Theta_t^{0,*},\Theta_t^{1,*},\mathbb{E}_t[\Theta_t^{1,*}];
L^0_t,L_t,P^0_t,Q^{00}_t,P_t,Q^{11}_t;u^1_t,\mathbb{E}_t[u^1_t]\big)
\big].$$
\end{proposition}

\subsection{Auxiliary control problem for the follower $\mathcal{A}_1$}
In this subsection, we fix $u^j, u^{0}$, the controls adopted by the leader $\mathcal{A}_j$ and the follower $\mathcal{A}_0$ respectively. Consider the control problem of the follower $\mathcal{A}_1$
associated with BSDE \eqref{081411-con}, SDE \eqref{mlf-state} and the functional \eqref{feb2902}. Given $(u^0,u^j),$ \eqref{081702-con} and the last two equations in \eqref{mlf-state} become exogenous processes.
The Hamiltonian functional $
H(t,x_0,y_0,z_0,u_0;x,y,z,u;\bar{x},\bar{y},\bar{z},\bar{u};l,p,q_{11})$
$=\langle p,b\rangle
+\langle q_{11},\sigma\rangle
- l\cdot f
+g,
$
and the related Hamiltonian system takes the following form
\begin{equation}\label{H-system-minor}
\left\{
\begin{aligned}
d{{X}}_t^{1,\ddagger}=&\partial_{p}H(t)dt
+\partial_{q_{11}}H(t)dW_t^1,\quad {X}_0^{1,\ddagger}=x_0,\\
-d{{Y}}_t^{1,\ddagger}=&-\partial_{l}H(t)dt
-{{Z}}_t^{1,\ddagger,0}dW_t^0-{{Z}}_t^{1,\ddagger,1}dW_t^1,\
Y_T^{1,\ddagger}=\Phi\big({{X}}_T^{1,\ddagger};
{X}_T^{0,\ddagger},\mathbb{E}_T[{{X}}^{1,\ddagger}_T]\big)+\xi^1,\\
dL_t^{\ddagger}=&-\partial_{y}H(t)dt
-\partial_{z}H(t)d(W^0,W^1)_t,\quad L_0^\ddagger=-\Gamma_y(Y_0^{1,\ddagger}),\\
-dP_t^\ddagger=&\partial_{x}H(t)dt-Q_t^{10,\ddagger}dW^0_t-Q_t^{11,\ddagger}dW_t,\quad
P_T^\ddagger=-L_T^\ddagger\cdot\partial_{{x}}\Phi({X}_T^{1,\ddagger};{X}_T^{0,\ddagger},
\mathbb{E}_T[{X}_T^{j,\ddagger}]),
\end{aligned}
\right.
\end{equation}
where the last two equations are the associated adjoint equations and
\begin{equation}\label{070602}
H(t)=\big(t,{\Theta_t^0},u_t^{0};{\Theta_t^1},u_t^{1,*};
{\mathbb{E}_t[\Theta_t^1]},\mathbb{E}_t[u_t^{1,*}];
L_t^\ddagger,P_t^\ddagger,Q^{11,\ddagger}_t\big).
\end{equation}Then from the stochastic maximum principle for FBSDE system (e.g. \cite{HLZ2023}), we have
\begin{proposition}\label{pro2}Under (A1), assume the FBSDE \eqref{H-system-minor} admits an unique solution and   \\
(1) There exists a unique maximizer of the Hamiltonian $H$ as a function of $u$
(denoted by $\widehat{u}$);\\
(2) The function $H$ is convex in $(x_0,x,y_0,y,z_0,z,u)$.\\
Then the optimal control of $\mathcal{A}_1$ is given by
$
\begin{aligned}
\widetilde{u}_t^1=\widehat{u}\Big(&t,\Theta_t^{0,\ddagger},
\Theta_t^{1,\ddagger},\mathbb{E}_t[\Theta^{j,\ddagger}_t];
L^{\ddagger}_t,P_t^\ddagger,Q^{11,\ddagger}_t;\mathbb{E}_t[u^j_t],u^0_t\Big).
\end{aligned}
$
\end{proposition}

\subsection{Consistency condition system}
For sake of presentation, hereafter we may write a bar on the top of a random variable (or process) to denote its conditional expectation with respect to $\mathbb{F}^0$, for example, $\overline{u}_t^j=\mathbb{E}_t[u_t^j]$.
We impose the following consistency conditions on the followers $\mathcal{A}_0, \mathcal{A}_1$ and the leader $\mathcal{A}_j$:
\begin{equation}\label{102901}
\widehat{u}\big[\overline{u}_t^j;\mathbb{E}_t[
\widehat{u}^{0}[u_t^1,\overline{u}_t^1]]\big]=u_t^1,\quad
\overline{u}^j_t=\overline{u}^{1}_t,\quad t\in[0,T].
\end{equation}
By \eqref{102901} and the solution uniqueness of SDE \eqref{mlf-state} and BSDEs \eqref{mlf-major-utility}, \eqref{081411-con}-\eqref{081702-con}, we have
\begin{equation}\label{011201}
\Theta^{0,\dagger}_t
=\Theta^{0,\ddagger}_t,\ \ \ \Theta^{1,\dagger}_t
=\Theta^{1,\ddagger}_t,\ \ \
\overline{\Theta}_t^{1,\dagger}
=\overline{\Theta}_t^{j,\ddagger},\ t\in[0,T].
\end{equation}
Noting \eqref{011201} and by Proposition \ref{pro1} and \ref{pro2}, we may introduce the following assumption:\\
\textbf{Assumption (A2)} \emph{
Suppose that there exist a pair of deterministic continuous functions
$({\varphi}^{0},{\varphi}):[0,T]\times\mathbb{R}^{3n+6m}\times \mathbb{R}^{3m}\times \mathbb{R}^{6n}\times U\rightarrow U^0\times U,$
satisfy the following conditions\begin{equation*}\label{011901}
\left\{
\begin{aligned}
\varphi^{0}=&
\widehat{u}^0\Big(t,\theta^{0},\theta,
\bar{\theta},
l^0,l,p^0,q^{00},p,q^{11},\bar{u};
\varphi\Big),\\
\varphi=&\widehat{u}\Big(t,\theta^{0},
\theta,
\bar{\theta}; l^{\ddagger},p^\ddagger,q^{11,\ddagger},
\bar{u};
\varphi^0\Big),\ t\in[0,T],
\end{aligned}
\right.
\end{equation*}
where $\widehat{u}^0$ and $\widehat{u}$ are the mappings given in Proposition \ref{pro1} and \ref{pro2}.
}

Under (A2), by measurable selection theorem,
there exists a measurable $\psi$ with
\begin{equation}\label{apr-14}
\psi_t=\mathbb{E}_t[\varphi(t,{\Theta_t^{0,\dagger},\Theta_t^{1,\dagger},
\overline{\Theta}^{1,\dagger}_t};
L^0_t,L_t,L_t^\ddagger,P^0_t,Q^{00}_t,P_t,Q^{11}_t,P_t^\ddagger,Q^{11,\ddagger}_t;\psi_t)].
\end{equation}
Combining \eqref{apr-14} and (A2), we can denote
\begin{equation}\label{apr-15}
\begin{aligned}
&\Psi^0(t,\Theta_t^{0,\dagger},\Theta_t^{1,\dagger};
L^0_t,L_t,L_t^\ddagger,P^0_t,Q^{00}_t,P_t,Q^{11}_t,P_t^\ddagger,Q^{11,\ddagger}_t)\\
=&
\mathbb{E}_t[\varphi^0(t,{\Theta_t^{0,\dagger},\Theta_t^{1,\dagger},
\overline{\Theta}^{1,\dagger}_t};
L^0_t,L_t,L_t^\ddagger,P^0_t,Q^{00}_t,P_t,Q^{11}_t,P_t^\ddagger,Q^{11,\ddagger}_t;\psi_t)],
\end{aligned}
\end{equation}
Similarly, $\Psi(t, \cdots)=
\varphi(t,{\Theta_t^{0,\dagger},\Theta_t^{1,\dagger},
\overline{\Theta}^{1,\dagger}_t};\cdots;\psi_t).$
Plugging the mappings $(\Psi^{0},\Psi)$ into \eqref{H-system-major} and \eqref{H-system-minor}, we obtain
\begin{equation}\label{ult-fbsde}
\left\{
\begin{aligned}
d{X}_t^{0}=&b^0\big(t,{X}_t^{0},
u_t^{0,*}; \mathbb{E}_t[{X}^{1}_t],
\mathbb{E}_t[{u}_t^{1,*}]\big)dt
+\sigma^0\big(t,{X}_t^{0},
u_t^{0,*}; \mathbb{E}_t[{X}^{1}_t],
\mathbb{E}_t[{u}_t^{1,*}]\big)dW_t^0,\\
d{{X}}_t^{1}=&b\big(t,X_t^{1},
u_t^{1,*};
{X}_t^{0},u_t^{0,*};
\mathbb{E}_t[X^{1}_t],
\mathbb{E}_t[u_t^{1,*}]\big)dt
+\sigma\big(t,X_t^{1},
u_t^{1,*};
{X}_t^{0},u_t^{0,*};
\mathbb{E}_t[X^{1}_t],
\mathbb{E}_t[u_t^{1,*}]\big)dW_t^1,\\
dL_t^0=&-\partial_{y_0}H_0(t)dt-\partial_{z_0}H_0(t)dW^0_t, \ \ \ \  \quad \quad  dL_t^{\ddagger}=-\partial_{y}H(t)dt
-\partial_{z}H(t)d(W^0,W^1)^\top_t,\\
dL_t=&-\big[\partial_{y}H_0(t)
+\mathbb{E}_t[\partial_{\bar{y}}H_0(t)]\big]dt
-\big[\partial_{z}H_0(t)+\mathbb{E}_t[\partial_{\bar{z}}H_0(t)]\big]
d(W^0,W^1)_t^\top,\\
-d{Y}_t^{0}=&f^0\big(t,{{\Theta}_t^{0},u_t^{0,*};
\mathbb{E}_t[{\Theta}^{1}_t],
\mathbb{E}_t[{{u}_t^{1,*}}]}\big)dt
-{Z}_t^{0}dW_t^0,\\
-d{{Y}}_t^{1}=&f\big(t,
{\Theta_t^{1},
{u_t^{1,*}};
\Theta_t^{0},u_t^{0,*};
\mathbb{E}_t[\Theta^{1}_t],\mathbb{E}_t[{u}_t^{1,*}]\big)dt
-{{Z}}_t^{1,0}dW_t^0-{{Z}}_t^{1,1}dW_t^1},\\
-dP_t^0=&\partial_{x_0}H_0(t)dt-Q_t^{00}dW^0_t-Q_t^{01}dW_t^1,\quad \quad -dP_t=\big[\partial_{x}H_0(t)
+\mathbb{E}_t[\partial_{\bar{x}}H_0(t)]\big]dt-Q_t^{10}dW^0_t-Q_t^{11}dW_t^1,\\
-dP_t^\ddagger=&\partial_{x}H(t)dt-Q_t^{10,\ddagger}dW^0_t
-Q_t^{11,\ddagger}dW_t^1,\\
\end{aligned}
\right.
\end{equation}with the mixed initial-terminal condition
$$
\begin{aligned}
{X}_0^{0}=&x_0^0,\ \ {X}_0^{1}=x_0,\ \ L_0^0=-\Gamma_y^0(Y_0^0),\ \ L_0^\ddagger=-\Gamma_y(Y_0^1),\ \ L_0=0, \\
{Y}_T^{0}=&
\Phi^0\big({X}_T^{0},
\mathbb{E}_T[{{X}}^{1}_T]\big)+\xi^0,\ \ \ \ {{Y}}_T^{1}=\Phi\big({{X}}_T^{1};
{X}_T^{0},{E}_T[{{X}}^{1}_T]\big)+\xi^1,\ \ \ \ P_T^\ddagger=-L_T^\ddagger\cdot\partial_{{x}}\Phi,\\
P_T^0=&-L_T^0\cdot\partial_{x_0}\Phi^0
-L_T\cdot\partial_{x_0}\Phi,\ \ \ \
P_T=-L_T\cdot\partial_x\Phi
-\mathbb{E}_T[L_T^0\cdot\partial_{\bar{x}}\Phi^0+L_T\cdot\partial_{\bar{x}}\Phi].
\end{aligned}
$$$H_0(t)$ (resp. $H(t)$) is given by \eqref{070601} (\eqref{070602}) by replacing $u^1$ ($u^0$) with $u^{1,*}$ ($u^{0,*}$) and
$$
\begin{aligned}
u_t^{0,*}=&\Psi^0(t,\Theta_t^{0,\dagger},\Theta_t^{1,\dagger};
L^0_t,L_t,L_t^\ddagger,P^0_t,Q^{00}_t,P_t,Q^{11}_t,P_t^\ddagger,Q^{11,\ddagger}_t),\\
u_t^{1,*}=&\Psi(t,\Theta_t^{0,\dagger},\Theta_t^{1,\dagger};
L^0_t,L_t,L_t^\ddagger,P^0_t,Q^{00}_t,P_t,Q^{11}_t,P_t^\ddagger,Q^{11,\ddagger}_t).
\end{aligned}
$$The following result is a direct consequence based on our previous analysis.
\begin{corollary}
Assume the assumptions of Proposition \ref{pro1} and \ref{pro2} hold. If FBSDE
\eqref{ult-fbsde} admits a solution, then $(u^{0,*},u^{1,*},u^{j,*})$ is an equilibrium of the triple-agent mixed leader-follower-Nash game, where
$u^{j,*}\in\mathcal{U}_j^d$  satisfies
$\overline{u}^{1,*}_t=\overline{u}^{j,*}_t$, for all $t\in[0,T]$.
\end{corollary}

\subsection{$\varepsilon$-Nash equilibrium of the RMM problem }
Proposition 4.3 yields an approximate Nash equilibrium for the RMM problem. To verify it, we first introduce an technical assumption, that is common in control literature \cite{CZ2016}:\\
\textbf{Assumption (A3)} \emph{
{\rm (i)} The diffusion coefficients $\sigma^0, \sigma$ are independent of $(u^0,u)$, if applicable.\\
{\rm (ii)} The maximizers $\widehat{u}^0$ and $\widehat{u}$ in Proposition \ref{pro1} and \ref{pro2} are independent of $(Z^0,Z^1)$.\\
{\rm (iii)} The  system \eqref{ult-fbsde} has
a unique solution $(\Theta_t^0,\Theta_t^1;L_t^0,L_t,L_t^\ddagger;P^0_t,P_t,$
$P^\ddagger_t,Q^{00}_t,Q_t^{01},Q_t^{10},Q_t^{11},$ $Q^{10,\ddagger}_t,
Q^{11,\ddagger}_t)$, where $L^0, L$ and $L^\ddagger$ are  $\mathbb{F}^0$-adapted.\\
\rm (iv) There exists a random decoupling field $\eta:[0,T]\times\Omega\times\mathbb{R}^{2n+3m}\rightarrow\mathbb{R}^{2m+3n}$
such that$$(Y^0_t,Y^1_t,P^0_t,P_t,P^\ddagger_t)=\eta(t,X_t^0,X_t^1,L_t^0,L_t,L_t^\ddagger),\ \text{a.s.},$$
with $\eta(\cdot,x^0,x^1,l^0,l,l^\ddagger)$ is $\mathbb{F}^0$-adapted for each $(x^0,x^1,l^0,l,l^\ddagger)\in\mathbb{R}^{2n+3m}$, and
$\eta(t, \cdots)$ is Lipschitz continuous on all its arguments uniformly on $t$.}

Decoupling field is a key to study the well-posedness of FBSDEs; see,
\cite{D2002,MWZZ2015,MY1999}. Proposition \ref{LQ-PRO} presents a sufficient condition ensuring the existence of decoupling field for LQG-RMM. By (A3-i), $\widehat{u}^0$ and $\widehat{u}$ in Proposition \ref{pro1} and \ref{pro2} are independent of $(Q^{00},Q^{11},Q^{11,\ddagger})$, and it is also the case for $\varphi^0,\varphi$ in (A2), $\psi$ in  \eqref{apr-14}, $\Psi^0$ and $\Psi$ in \eqref{apr-15}. We abuse the notation to write
\begin{equation}\nonumber
\begin{aligned}
&\chi(t,X^{0}_t,
X^{1}_t)=\chi(t,X^{0}_t,
X^{1}_t;L_t^0,L_t,
L_t^\ddagger,
\eta(t,X^{0}_t,
X^{1}_t,L_t^0,L_t,
L_t^\ddagger)),\ \chi=\Psi^0, \Psi.\\
\end{aligned}
\end{equation}
Then it follows from \eqref{apr-15} and (A3-iii) that $\Psi^0(t,X_t^0,X_t^1)$ and $\Psi(t,\cdot,\cdot)$ are $\mathbb{F}^0$-adapted.

\textbf{Assumption (A4)} \emph{Both $\Psi^0$ and $\Psi$ are Lipschitz in state variable
$(x^0,x)$ uniformly on $t$, and $||\Psi^0(t,0,0)||_{L^{2}}+||\Psi(t,0,0)||_{L^{2}}<\infty.$}

(A3) and (A4) are commonly adopted in MFG literature (e.g., (A7) and (A8) in \cite{CZ2016}). Indeed, (A3) and (A4) are satisfied for the LQG-RMM problem in Section 5.
Applying the feedback control pair $(\Psi^0,\Psi)$ to \eqref{major-state} and \eqref{major-utility}, we get the following forward-backward system
\begin{equation}\label{optimal-state}
\left\{
\begin{aligned}
dX_t^{0,N}=&b^0\big(t,X_t^{0,N},
u_t^{0,N}; X_t^{(N)},{u_t^{(N)}}\big)dt+{\sigma^{0}\big(t,X_t^{0,N}; X_t^{(N)},{u_t^{(N)}}\big)}dW^0_t,\
X_0^0=x_0^0,\\
dX_t^{i,N}=&b\big(t,X_t^{i,N},
u_t^{i,N};
X_t^{0,N},u_t^{0,N};X_t^{(N)},u_t^{(N)}\big)dt\\
&+\sigma\big(t,X_t^{i,N},X_t^{0,N};X_t^{(N)},u_t^{(N)}\big)dW_t^i,\quad
X_0^i=x_0,\\
-dY_t^{0,N}=&f^0\big(t,\Theta_t^{0,N},
u_t^{0,N};  \Theta_t^{(N)},u_t^{(N)}\big)dt-Z_t^{0,N}dW_t^0,\
Y_T^{0,N}=\Phi^0(X_T^{0,N},X_T^{(N)})+\xi^0,\\
-dY_t^{i,N}=&f\big(t,\Theta_t^{i,N},
u_t^{i,N};
\Theta_t^{0,N},
u_t^{0,N};
 \Theta_t^{(N)},u_t^{(N)}\big)dt
-Z_t^{i,N}d(W^0,W^i)_t,\\
Y_T^i=&\Phi(X_T^{i,N},X_T^{0,N},X_T^{(N)})+\xi^i,
\end{aligned}
\right.
\end{equation}
where $\Xi_\cdot^{(N)}=\frac1N\sum_{i=1}^N\Xi_\cdot^{i,N},\ \Xi=X,Y,Z,\  u_\cdot^{(N)}=\frac1N\sum_{i=1}^Nu_t^{i,N}$,
$u_t^{0,N}=\Psi^0(t,X_t^{0,N},
X_t^{1,N}),$ $u_t^{i,N}=
\Psi(t,X_t^{0,N};X_t^{i,N}).$
The payoff functionals \eqref{091601} now become
\begin{equation}\label{optimal-value}
\begin{aligned}
J_0^N=\Gamma^0(Y_0^{0,N})+\mathbb{E}\big[
\int_0^Tg^0&\big(t,\Theta_t^{0,N}, u_t^{0,N};\Theta_t^{(N)},u_t^{(N)}\big)dt\big],\\
J_{i}^N
=\Gamma(Y_0^{i,N})+\mathbb{E}\big[\int_0^T
g&\big(t,\Theta_t^{i,N},u_t^{i,N};
\Theta_t^{0,N}, u_t^{0,N};\Theta_t^{(N)},u_t^{(N)}\big)dt\big].
\end{aligned}
\end{equation}To ease the notation, hereafter we may write $(J_0^N,J_i^N)$ instead of $(J_0^N,J_i^N)(u^0,u^1,\cdots,u^N)$.

\begin{theorem}\label{th4.1}
Assume that (A1)-(A4) hold. The feedback control
$\big(
u^{0,N}_t,
u^{1,N}_t,\cdots,u^{N,N}_t\big)$
is an $\varepsilon_N$-Nash equilibrium for RMM problem, where $\varepsilon_N\leq \frac{C}{\sqrt{N}}$.
\end{theorem}
The proof of Theorem \ref{th4.1} essentially relies on a tailor-made propagation chaos analysis in the same spirit of \cite{CZ2016}, but invokes additional estimates on BSDEs. We defer it in Appendix \ref{app1}.

\section{Linear-quadratic-Gaussian cases}\label{LQ}
{This section studies the RMM problem in some linear-quadratic-Gaussian (LQG) cases (LQG-RMM), where the agents' states evolve by the following linear SDEs: for $1\leq i\leq N$,
\begin{equation}\label{LQ-state}
\left\{
\begin{aligned}
dX_t^0=&\big(b_1^0X_t^0+b_2^0{u_t^0}+b_3^0X_t^{(N)}+b_4^0{u_t^{(N)}}\big)dt
+\sigma_0dW^0_t,\quad X_0^0=x_0^0,\\
dX_t^i=&\big(b_1X_t^i+b_2u_t^i+
b_3X_t^0+b_4u_t^0+b_5X_t^{(N)}+b_6u_t^{(N)}\big)dt
+\sigma dW_t^i,\quad X_0^i=x_0,
\end{aligned}
\right.
\end{equation}and the recursive functionals assume the following quadratic forms
\begin{equation}\label{feb2903}
\begin{aligned}
&J_0^N
=-\gamma_0|Y_0^0|^2-\mathbb{E}\big[\int_0^T Q_0|X_t^0-\mu^0_1\cdot X_t^{(N)}|^2+ R_0|u_t^0-\mu^0_2\cdot u_t^{(N)}|^2dt\big],\\
&J_{i}^N
=-\gamma|Y_0^i|^2
-\mathbb{E}\big[\int_0^T Q|X_t^i-\mu_1\cdot X_t^{(N)}-\mu_2\cdot X_t^{0}|^2
+ R|u_t^i-\mu_3\cdot u_t^{(N)}-\mu_4\cdot u_t^0|^2dt\big],
\end{aligned}
\end{equation}where $(Y^0,Z^0;Y^i,Z^i)$ satisfy the following coupled linear BSDE system:
\begin{equation}\label{LQ-payoff}
\left\{
\begin{aligned}
-dY_t^0=&\big( f_1^0X_t^0+f_2^0Y_t^0+f_3^0Z_t^0
+ f_4^0u_t^0
+f_5^0X_t^{(N)}+f_6^0Y_t^{(N)}+f_7^0Z_t^{(N),0}+f_8^0 u_t^{(N)}\big)dt-Z_t^0dW_t^0,\\
-dY_t^i=&\big( f_1X_t^i+f_2 Y_t^{i}+f_3 Z_t^{i,0}
+f_4u_t^i
+f_5X_t^0+f_6 Y_t^{0}+f_7 Z_t^{0}+f_8 u_t^{0}
+f_9X_t^{(N)}\\
&+f_{10}Y_t^{(N)}+f_{11}Z_t^{(N),0}+f_{12} u_t^{(N)}
\big)dt
-Z_t^{i,0}dW_t^0-Z_t^{i,i}dW_t^i,\\
Y_T^0=&\Phi_1^0X_T^0+\Phi_2^0 X_T^{(N)}+\xi^0,\quad
Y_T^i= \Phi_1X_T^i +\Phi_2X_T^0 +\Phi_3X_T^{(N)}+\xi^i.
\end{aligned}
\right.
\end{equation}with $ Y_t^{(N)}=\frac{1}{N}\sum_{i=1}^NY_t^{i}, Z_t^{(N),0}=\frac{1}{N}\sum_{i=1}^NZ_t^{i,0}$. To simplify the analysis, we assume that $n=m=k^0=k=1$; all coefficients are constants with nonnegative $(\gamma_0,\gamma, Q_0,Q)$, positive $(R_0,R)$, and {$\mu_3+\mu^0_2\mu_4\neq 1$}.
The forward-backward LQG setting of \eqref{LQ-state}-\eqref{LQ-payoff} is strongly motivated by various practical applications, see \cite{HJX2023, KEY2020, MW1995} for more details. We introduce the following notations:
\begin{equation}\label{052301}
\begin{aligned}
&\mathbf{X}_t=(X^0_t,{X}^1_t)^\top,\
\mathbf{L}_t=({L}_t^0,{L}_t,{L}^\ddagger_t)^\top,\
\mathbf{Y}_t=(Y_t^0,{Y}_t^1)^\top,\
\mathbf{P}_t=({P}^0_t,{P}_t,{P}^\ddagger_t)^\top,\\
&\mathbf{Z}_t=(Z_t^0,{Z}_t^{1,0})^\top,\
\mathbf{Q}_t=({Q}^{00}_t,{Q}^{10}_t,{Q}^{10,\ddagger}_t)^\top,\
\mathbf{Z}_t^1=(0,{Z}_t^{1,1})^\top,\
\mathbf{Q}_t^1=({Q}^{01}_t,{Q}^{11}_t,{Q}^{11,\ddagger}_t)^\top.
\end{aligned}
\end{equation}By Proposition \ref{pro1}, \ref{pro2} and \eqref{102901}, the equilibrium strategy $(u^{0,*},u^{1,*},u^{j,*})$ now read as
\begin{equation}\label{11150101}
\begin{aligned}
 u_t^{0,*}=&{a_1^0}^\top \mathbb{E}_t[\mathbf{P}_t]+{a_3^0}^\top \mathbb{E}_t[\mathbf{L}_t],\quad
u_t^{1,*}=
a_1^\top \mathbb{E}_t[\mathbf{P}_t]+a_3^\top \mathbb{E}_t[\mathbf{L}_t]+a_7^\top \mathbf{P}_t+a_8^\top \mathbf{L}_t,\\
 \mathbb{E}_t[{u}_t^{j,*}]=&({a}_1+a_7)^\top \mathbb{E}_t[\mathbf{P}_t]+({a}_3+a_8)^\top \mathbb{E}_t[\mathbf{L}_t].
 \end{aligned}
 \end{equation}Combining with \eqref{ult-fbsde}, we get the following FBSDE
\begin{equation}\label{043001}
\left\{
\begin{aligned}
d\mathbf{X}_t=&
\Big(\mathbb{A}_1\mathbf{X}_t
+\mathbb{A}_2\mathbb{E}_t[\mathbf{X}_t]
+\mathbb{B}_1\mathbf{P}_t
+\mathbb{B}_2\mathbb{E}_t[\mathbf{P}_t]
+\mathbb{C}_1\mathbf{L}_t
+\mathbb{C}_2\mathbb{E}_t[\mathbf{L}_t]
\Big)dt\\
&+(\sigma_0,0)^\top dW_t^0+(0,\sigma)^\top dW_t^1,\
\mathbf{X}_0=({x}_0^0,x_0)^\top,\\
d\mathbf{L}_t=&
\Big(\mathbb{C}_3\mathbf{L}_t
+\mathbb{C}_4\mathbb{E}_t[\mathbf{L}_t]
\Big)dt+\Big(\mathbb{C}_5\mathbf{L}_t+\mathbb{C}_6\mathbb{E}_t[\mathbf{L}_t]
\Big)dW_t^0,\ \mathbf{L}_0=\rho\mathbf{Y}_0,\\
-d\mathbf{Y}_t=&\Big(\mathbb{A}_3\mathbf{X}_t
+\mathbb{A}_4\mathbb{E}_t[\mathbf{X}_t]
+\mathbb{D}_1\mathbf{Y}_t+\mathbb{D}_2\mathbb{E}_t[\mathbf{Y}_t]
+\mathbb{F}_1\mathbf{Z}_t
+\mathbb{F}_2\mathbb{E}_t[\mathbf{Z}_t]\\
&+\mathbb{B}_3\mathbf{P}_t
+\mathbb{B}_4\mathbb{E}_t[\mathbf{P}_t]
+\mathbb{C}_7\mathbf{L}_t
+\mathbb{C}_8\mathbb{E}_t[\mathbf{L}_t]
\Big)dt
-\mathbf{Z}_tdW_t^0-\mathbf{Z}_t^1dW_t^1,\\
-d\mathbf{P}_t=&\Big(\mathbb{A}_5\mathbf{X}_t
+\mathbb{A}_6\mathbb{E}_t[\mathbf{X}_t]
+\mathbb{B}_5\mathbf{P}_t+\mathbb{B}_6\mathbb{E}_t[\mathbf{P}_t]
+\mathbb{C}_9\mathbf{L}_t
+\mathbb{C}_{10}\mathbb{E}_t[\mathbf{L}_t]
\Big)dt
-\mathbf{Q}_tdW_t^0-\mathbf{Q}_t^1dW_t^1,\\
\mathbf{Y}_T=&\mathbb{A}_7\mathbf{X}_T+\mathbb{A}_8\mathbb{E}_T[\mathbf{X}_T]+(\xi^0,\xi)^\top,\
\mathbf{P}_T=\mathbb{C}_{11}\mathbf{L}_T+\mathbb{C}_{12}\mathbb{E}_T[\mathbf{L}_T].
\end{aligned}
\right.
\end{equation}For sake of presentation, we defer the definitions of $\{a_i^0\}_{i=1,3}$, $\{a_i\}_{i=1,3,7,8}$ and matrices (vectors) of \eqref{11150101}-\eqref{043001} to Appendix \ref{app2}.  \eqref{043001} is a fully-coupled FBSDE with mixed initial-terminal conditions, and its wellposedness can be discussed through the following steps. Set
$
\overline{M}_t
=\mathbb{E}_t[M_t]$ for $M=\mathbf{X},\mathbf{Y},\mathbf{Z},\mathbf{L},\mathbf{P},\mathbf{Q},
$
and take conditional expectation on \eqref{043001}, we have
\begin{equation}\label{111405}
\left\{
\begin{aligned}
d\overline{\mathbf{X}}_t=&[(\mathbb{A}_1+\mathbb{A}_2)\overline{\mathbf{X}}_t
+(\mathbb{B}_1+\mathbb{B}_2)\overline{\mathbf{P}}_t
+(\mathbb{C}_1+\mathbb{C}_2)\overline{\mathbf{L}}_t]dt+(\sigma_0,0)^\top dW_t^0,\quad
\overline{\mathbf{X}}_0=({x}_0^0,x_0)^\top,\\
d\overline{\mathbf{L}}_t=&
(\mathbb{C}_3
+\mathbb{C}_4)\overline{\mathbf{L}}_tdt
+(\mathbb{C}_5+\mathbb{C}_6)\overline{\mathbf{L}}_tdW_t^0,\quad
\overline{\mathbf{L}}_0=\rho\overline{\mathbf{Y}}_0,\\
-d\overline{\mathbf{Y}}_t=&[(\mathbb{A}_3+\mathbb{A}_4)\overline{\mathbf{X}}_t
+(\mathbb{D}_1+\mathbb{D}_2)\overline{\mathbf{Y}}_t
+(\mathbb{F}_1+\mathbb{F}_2)\overline{\mathbf{Z}}_t
+(\mathbb{B}_3+\mathbb{B}_4)\overline{\mathbf{P}}_t
+(\mathbb{C}_7+\mathbb{C}_8)\overline{\mathbf{L}}_t]dt
-\overline{\mathbf{Z}}_tdW_t^0,\\
-d\overline{\mathbf{P}}_t=&[(\mathbb{A}_5
+\mathbb{A}_6)\overline{\mathbf{X}}_t
+(\mathbb{B}_5+\mathbb{B}_6)\overline{\mathbf{P}}_t
+(\mathbb{C}_9
+\mathbb{C}_{10})\overline{\mathbf{L}}_t]dt
-\overline{\mathbf{Q}}_tdW_t^0,\\
\overline{\mathbf{Y}}_T=&(\mathbb{A}_7+\mathbb{A}_8)\overline{\mathbf{X}}_T
+(\xi^0,\mathbb{E}_T[\xi])^\top,\quad
\overline{\mathbf{P}}_T=(\mathbb{C}_{11}+\mathbb{C}_{12})\overline{\mathbf{L}}_T.
\end{aligned}
\right.
\end{equation}
Next, we consider the following ODE and BSDE
\begin{equation}\label{31204}
\left\{
\begin{aligned}
&\dot{S}_t+{S}_t{\mathbb{H}_1}
+\mathbb{H}_2S_t+S_t{\mathbb{H}_3}S_t
+
(I+S_t\widetilde{\rho})\mathbb{H}_4(I+S_t\widetilde{\rho})^{-1}
\big(S_t\mathbb{H}_5
+S_t\mathbb{H}_6 S_t\big)+
\mathbb{H}_7=0,\\
&{S}_T=(I-\mathbb{G}_1)^{-1}\mathbb{G}_2,
\end{aligned}
\right.
\end{equation}
\begin{equation}\label{31205}
\left\{
\begin{aligned}
-d\Upsilon_t=&\Big[\Big(S_t \mathbb{H}_3+\mathbb{H}_2
+(I+S_t\widetilde{\rho})\mathbb{H}_4(I+S_t\widetilde{\rho})^{-1}S_t\mathbb{H}_6\Big)\Upsilon_t
+(I+S_t\widetilde{\rho})\mathbb{H}_4(I+S_t\widetilde{\rho})^{-1}\upsilon_t\\
&+(S_t+I)\mathbb{H}_4(I+S_t\widetilde{\rho})^{-1}S_t(\sigma_0,0,0,0,0)^\top\Big]dt-\upsilon_tdW_t^0,\\
\Upsilon_T=&(I-\mathbb{G}_1)^{-1}(\xi^0,\mathbb{E}_T[\xi],0,0,0)^\top.
\end{aligned}
\right.
\end{equation}
Again, we defer the definitions of $\{\mathbb{H}_k\}_{k=1}^7$ and $\mathbb{G}_1,\mathbb{G}_2$ to Appendix \ref{app2}. The well-posedness of \eqref{31204} may be obtained as Theorem 4.6 in   \cite{HJX2023}. We refrain to present these conditions in details, as they might be rather technical and incur unnecessary degression along our presentation. Instead, we directly assume that

(A5) \eqref{31204} admits a unique $\mathbb{R}^{5\times 5}$-valued solution $S$ with bounded $[I+S_t\widetilde{\rho}]^{-1},\  \text{a.e.}\ t.$

Under (A5), BSDE \eqref{31205} admits an unique $\mathbb{F}^0$-adapted solution $(\Upsilon,\upsilon)$ by the standard BSDE solvability arguments. Moreover, we have the following result:\begin{lemma}\label{le5.1}
Let (A5) holds. Then the linear FBSDE \eqref{111405} has a unique $\mathbb{F}^0$-adapted solution $(\overline{\mathbf{X}},\overline{\mathbf{L}},\overline{\mathbf{Y}},\overline{\mathbf{Z}},
\overline{\mathbf{P}},\overline{\mathbf{Q}})$ with the following relations:\begin{equation}\label{31203}
\begin{aligned}
(\overline{\mathbf{Y}}_t,\overline{\mathbf{P}}_t)^\top=&
(I+S_t\widetilde{\rho})^{-1}S_t(\overline{\mathbf{X}}_t,\overline{\mathbf{L}}_t)^\top
+(I+S_t\widetilde{\rho})^{-1}\Upsilon_t,\\ (\overline{\mathbf{Z}}_t,\overline{\mathbf{Q}}_t)^\top=&
(I+S_t\widetilde{\rho})^{-1}(S_t\mathbb{H}_5+S_t
\mathbb{H}_6S_t)
\big(I-\widetilde{\rho}(I+S_t\widetilde{\rho})^{-1}S_t\big)
(\overline{\mathbf{X}}_t,\overline{\mathbf{L}}_t)^\top\\
&-(I+S_t\widetilde{\rho})^{-1}(S_t\mathbb{H}_5
+S_t\mathbb{H}_6S_t)\widetilde{\rho}(I+S_t\widetilde{\rho})^{-1}\Upsilon_t\\
&+(I+S_t\widetilde{\rho})^{-1}S_t\mathbb{H}_6\Upsilon_t+(I+S_t\widetilde{\rho})^{-1}
(S_t(\sigma_0,0,0,0,0)^\top+\upsilon_t).
\end{aligned}
\end{equation}
\end{lemma}Its proof is based on a standard linear transformation decoupling method ({e.g., \cite{FHH2022}) and we omit its details here. We now assume

(A6)
$
 \det\{(0,I)e^{
 \mathcal{A}t}(0,I)^\top\}>0,\ \forall\ t\in[0,T],
$
where $\mathcal{A}=\left(
 \begin{array}{cc}
 b_1&\frac12 R^{-1}|b_2|^2\\
 2Q&-b_1
 \end{array}
 \right)$. \begin{proposition}\label{LQ-PRO}
Let (A5)-(A6) hold. Then the Hamiltonian system \eqref{043001} has a unique solution
 $(\mathbf{X},\mathbf{L},\mathbf{Y},\mathbf{Z},\mathbf{Z}^1,\mathbf{P},\mathbf{Q},\mathbf{Q}^1)$ (see \eqref{052301} for notation), where $(X^0,\mathbf{L},Y^0,Z^0,P^0,Q^{00},P,Q^{10})$ are $\mathbb{F}^0$-adapted, $(X^1,Y^1,Z^{10},Z^{11},P^\ddagger,Q^{10,\ddagger},Q^{11,\ddagger})$ are $\mathbb{F}^{0,1}$-adapted, and $Q^{01}=Q^{11}\equiv 0$. Moreover, the following relation holds
 $
P_t^\ddagger=\Sigma_t X_t^{1}+p_t,\ t\in[0,T]$, where $\Sigma$ and $p$ are given by
\begin{equation}\label{042503}
\Sigma_t=-\Big[(0,I)e^{
 \mathcal{A}(T-t)}(0,I)^\top\Big]^{-1}(0,I)e^{
 \mathcal{A}(T-t)}(I,0)^\top,
\end{equation}
\begin{equation}\label{042601}
\begin{aligned}
&p_t=\mathbb{E}_t\Big[-L_T^\ddagger\Phi_1\Pi_T
+
\int_t^T[(\Sigma_s\mathbf{\Lambda}_1
+\Sigma_s\mathbf{\Lambda}_2(I+S_s\widetilde{\rho})^{-1}S_s+\mathbf{\Lambda}_4)(\overline{\mathbf{X}}_s,
\mathbf{L}_s)^\top+
\Sigma_s\mathbf{\Lambda}_2(I+S_s\widetilde{\rho})^{-1}\Upsilon_s]
\Pi_sds\Big],
\end{aligned}
\end{equation}with $\Pi_s=e^{\int_t^s(\frac12 R^{-1}|b_2|^2\Sigma_r+b_1)dr}, \ \ s\in[t,T]$, and $\mathbf{\Lambda}_1,\mathbf{\Lambda}_2$ and $\mathbf{\Lambda}_4$ are given in Appendix \ref{app2}.
\end{proposition}
We defer the proof of Proposition \ref{LQ-PRO} in Appendix. Combining \eqref{11150101}, \eqref{31203} and \eqref{042601}, the equilibrium for  $\mathcal{A}_0$ and $\mathcal{A}_1$ are respectively
\begin{equation*}\label{042505}
\begin{aligned}
u_t^{0,*}=&(0,{a_1^0}^\top) (\overline{\mathbf{Y}}_t,\overline{\mathbf{P}}_t)^\top+{a_3^0}^\top {\mathbf{L}}_t=A_t^{0,1}(X_t^0,\mathbb{E}_t[X_t^1])^\top+M_t^0,\\
u_t^{1,*}=&
(0,a_1^\top)(\overline{\mathbf{Y}}_t,\overline{\mathbf{P}}_t)^\top+(a_3+a_8)^\top \mathbf{L}_t+\frac12 R^{-1}b_2 P_t^\ddagger
=A_t^{1}(X_t^0,\mathbb{E}_t[X_t^1])^\top\\
&+\frac12 R^{-1}b_2\Sigma_tX_t^1
+\frac12 R^{-1}b_2\mathbb{E}_t\Big[
\int_t^TC_s^{1}(X_s^0,\mathbb{E}_s[X_s^1])^\top]
\Pi_sds\Big]+M_t,
\end{aligned}
\end{equation*}
where
\begin{equation*}\label{042602}
\begin{aligned}
&(A_t^{0,1},A_t^{0,2}):=
(0,{a_1^0}^\top)(I+S_t\widetilde{\rho})^{-1}
S_t,\ \ \ \ \ M_t^0:=A_t^{0,2}\mathbf{L}_t+(0,{a_1^0}^\top)(I+S_t\widetilde{\rho})^{-1}\Upsilon_t+{a_3^0}^\top \overline{\mathbf{L}}_t,\\
&(A_t^{1},A_t^{2}):=(0,a_1^\top)(I+S_t\widetilde{\rho})^{-1}
S_t,\ \ \ \ \ \
(C_t^{1},C_t^{2})
:=\Sigma_t\mathbf{\Lambda}_1+\Sigma_t\mathbf{\Lambda}_2
(I+S_t\widetilde{\rho})^{-1}S_t+\mathbf{\Lambda}_4,\\
\end{aligned}
\end{equation*}
\begin{equation*}\label{042602}
\begin{aligned}
&M_t:=A_t^{2}\mathbf{L}_t
+(0,a_1^\top)(I+S_t\widetilde{\rho})^{-1}\Upsilon_t+(a_3+a_8)^\top\mathbf{L}_t\\
&\qquad+\frac12 R^{-1}b_2\mathbb{E}_t\Big[-L_T^\ddagger\Phi_1\Pi_T
+
\int_t^T[C_s^2\mathbf{L}_s+
\Sigma_s\mathbf{\Lambda}_2(I+S_s\widetilde{\rho})^{-1}\Upsilon_s]
\Pi_sds\Big].
\end{aligned}
\end{equation*}
Note that both $N_t^0$ and $N_t$ are $\mathbb{F}^0$-adapted.
By Theorem \ref{th4.1}, we have the following result.
\begin{theorem}\label{th5.1}
Let (A5)-(A6) hold. Then
\begin{equation}\label{0915}
\begin{aligned}
\Big(&
A_t^{0,1}(X_t^{0,N},
\mathbb{E}_t[X_t^{1,N}])^\top+M_t^0,\ \ \
A_t^{1}(X_t^{0,N},\mathbb{E}_t[X_t^{1,N}])^\top\\
&+\frac12 R^{-1}b_2\Sigma_tX_t^{i,N}
+\frac12 R^{-1}b_2\mathbb{E}_t\Big[
\int_t^TC_s^{1}(X_s^{0,N},\mathbb{E}_s[X_s^{1,N}])^\top]
\Pi_sds\Big]+M_t\Big),\ 1\leq i\leq N,
\end{aligned}
\end{equation}is an $\varepsilon_N$-Nash equilibrium {strategy} for LQG-RMM with $\varepsilon_N= O(\frac{1}{\sqrt{N}})$; $(X^{0,N},X^{1,N},$ $\cdots,X^{N,N})$ is the solution of McKean-Vlasov SDE \eqref{LQ-state} by applying the equilibrium \eqref{0915}.
\end{theorem}

 \subsection{Forward LQG-RMM}
This subsection studies a special case with $\gamma_0=\gamma=0$ in \eqref{LQ-state}, \eqref{feb2903}. In this case, all functionals are still quadratic but involve only the \emph{forward} state, and the RMM reduces to the classical \emph{forward} major-minor game in \cite{CZ2016}}. Although such forward setting is not novel, our RMM-LQG still gains novelties in its generality of weak-couplings. Now \eqref{31204} and \eqref{31205} become\begin{equation*}\label{042604}
\left\{
\begin{aligned}
&\dot{S}_t+{S}_t(\mathbb{A}_1+\mathbb{A}_2)+(\mathbb{B}_5+\mathbb{B}_6)S_t
+S_t\mathbb{B}_2S_t+(\mathbb{A}_5+\mathbb{A}_6)=0,\\
&-d\Upsilon_t=(S_t\mathbb{B}_2+
\mathbb{B}_5+\mathbb{B}_6)\Upsilon_tdt,\ \ \ \
S_T=\Upsilon_T=0.
\end{aligned}
\right.
\end{equation*}It is easy to check that if\begin{equation}\label{A3}
 \det\{(0,I)e^{
 \widehat{\mathcal{A}}t}(0,I)^\top\}>0,\ \forall\ t\in[0,T],\ \ \ \text{with} \  \widehat{\mathcal{A}}=\left(
 \begin{array}{cc}
 \mathbb{A}_1+\mathbb{A}_2&\mathbb{B}_2\\
 -\mathbb{A}_5-\mathbb{A}_6&-\mathbb{B}_5-\mathbb{B}_6
 \end{array}
 \right),
 \end{equation}
 then (A5) can be ensured and $S$ admits the following representation\begin{equation}\label{042605}
S_t=-\Big[(0,I)e^{
 \widehat{\mathcal{A}}(T-t)}(0,I)^\top\Big]^{-1}(0,I)e^{
 \widehat{\mathcal{A}}(T-t)}(I,0)^\top=:
 \left(
 \begin{array}{ccc}
 S_t^{11}& S_t^{12}\\
  S_t^{21}& S_t^{22}\\
   S_t^{31}& S_t^{32}\\
 \end{array}
 \right).
\end{equation}By the solution uniqueness, $
\Upsilon_t\equiv 0,
$ so $M_t^0=M_t\equiv 0, \ t\in[0,T]$ (noting $L^0=L=L^\ddagger\equiv 0$).
Then
\begin{corollary}\label{le5.2}
Let  (A6) and  \eqref{A3} hold. An $\varepsilon_N$-Nash equilibrium {strategy} for the forward LQG-RMM problem takes the following form$$
\begin{aligned}
\Big(&
A_t^{0,1}(X_t^{0,N},
\mathbb{E}_t[X_t^{1,N}])^\top,\
A_t^{1}(X_t^{0,N},\mathbb{E}_t[X_t^{1,N}])^\top\\
&+\frac12 R^{-1}b_2\Sigma_tX_t^{i,N}
+\frac12 R^{-1}b_2\mathbb{E}_t\Big[
\int_t^TC_s^{1}(X_s^{0,N},\mathbb{E}_s[X_s^{1,N}])^\top]
\Pi_sds\Big]\Big),\quad 1\leq i\leq N,
\end{aligned}
$$where
$A_t^{0,1}={a_1^0}^\top (S_t^{11},S_t^{21},S_t^{31})^\top,\
A_t^{1}=a_1^\top (S_t^{11},S_t^{21},S_t^{31})^\top,\ C_t^{1}=\Sigma_t(b_3,b_5)+\Sigma_t\mathbf{\Lambda}_3S_t+(2Q \mu_2,2Q\mu_1),$
 $\mathbf{\Lambda}_3$ is given in Subsection 6.2; $S_t$ and $\Sigma_t$ are given by \eqref{042605} and \eqref{042503}, respectively.
\end{corollary}
\begin{example}\label{eg5.1}Consider a forward LQG-RMM with $\mu^0_2=\mu_3=\mu_4=0$ and $b_4^0=b_4=b_6=0$ (no weak-coupling of the control-average). The $\varepsilon_N$-Nash equilibrium by Corollary 5.1 becomes
$$
\begin{aligned}
\Big(
\frac12 R_0^{-1}b_2^0S_t^{11}X_t^{0,N}+\frac12 R_0^{-1}b_2^0S_t^{21}\mathbb{E}_t[X_t^{1,N}],\ \ \
\frac12 R^{-1}b_2 \Sigma_t X_t^{i,N}+k_t\Big),\ 1\leq i\leq N,
\end{aligned}
$$
where $k_t=\frac12 R^{-1}b_2\mathbb{E}_t\Big[
\int_t^TC_s^{1}(X_s^{0,N},\mathbb{E}_s[X_s^{1,N}])^\top\Pi_sds\Big]$ and $C_t^{1}=\Sigma_t(b_3,b_5)+(2Q \mu_2,2Q\mu_1).$
This result recovers Theorem 5.1 in \cite{CZ2016} with a subtle difference: note that $\varepsilon_N=O(N^{-\frac{1}{d+4}})$ in \cite{CZ2016} while our $\varepsilon_N=O(\sqrt{N})$. This is mainly due to the modeling differences in weak-couplings: \cite{CZ2016} considers the empirical distribution while we focus on the more special empirical average. As a tradeoff, we can obtain an  explicit expression for $k_t$ while \cite{CZ2016} only shows its existence.
\end{example}

\subsection{Backward LQG-RMM} In this subsection, $Q_0=Q=0$, $f_1^0=f_5^0=f_{1,5,9}=0$, $\Phi_1^0=\Phi_2^0=\Phi_{1,2,3}=0$ in \eqref{feb2903} and \eqref{LQ-payoff}. The LQG-RMM is now solely ``backward" without the forward state \eqref{LQ-state}.
Backward LQG setting has found broad applications in such as optimal investment, recursive utility and hedging (e.g., \cite{HHL2022}). Related MFG studies on this setting has also been well addressed (see \cite{DHW2018,HWW2016}). However, the backward LQG-RMM seems still novel in literature. It can well capture the insights of the \emph{large investor} (see \cite{EPQ2001}) and \emph{relative performance} (\cite{DHW2018}), both are well motivated in financial studies. In this case, \eqref{31204} reads as (noting that $S$ is $\mathbb{R}^{2\times 3}$-valued)
\begin{equation}\label{042606}
\begin{aligned}
&\dot{S}_t+{S}_t(\mathbb{C}_3+\mathbb{C}_4+\rho\mathbb{C}_7+\rho\mathbb{C}_8)
+\widetilde{\mathbb{C}}S_t
+S_t\widetilde{\mathbb{B}}S_t+
(I+S_t\rho)(\mathbb{F}_1+\mathbb{F}_2)(I+S_t\rho)^{-1}\cdot\\
&\big(S_t(\mathbb{C}_5+\mathbb{C}_6)+S_t(\mathbb{C}_5+\mathbb{C}_6)\rho S_t\big)+
\mathbb{C}_7+\mathbb{C}_8=0,\
{S}_T=0,
\end{aligned}
\end{equation}and the linear BSDE \eqref{31205} now becomes
\begin{equation}\label{042607}
\left\{
\begin{aligned}
-d\Upsilon_t=&\Big[\Big(S_t\widetilde{\mathbb{B}}+
\widetilde{\mathbb{C}}
+(I+S_t\rho)(\mathbb{F}_1+\mathbb{F}_2)(I+S_t\rho)^{-1}S_t
(\mathbb{C}_5+\mathbb{C}_6)\rho\Big)\Upsilon_t\\
&+
(I+S_t\rho)(\mathbb{F}_1+\mathbb{F}_2)(I+S_t\rho)^{-1}\upsilon_t
\Big]dt-\upsilon_tdW_t^0,\
\Upsilon_T=(\xi^0,\mathbb{E}_T[\xi^1])^\top,
\end{aligned}
\right.
\end{equation}
where  $\widetilde{\mathbb{C}}=
\mathbb{C}_7\rho+\mathbb{C}_8\rho+\mathbb{D}_1+\mathbb{D}_2,$
$\widetilde{\mathbb{B}}=(\mathbb{C}_3+\mathbb{C}_4)\rho
+\rho\widetilde{\mathbb{C}}$.
Moreover, it follows from \eqref{042503}-\eqref{042601} that $\Sigma_t=p_t\equiv 0$.
Besides, by \eqref{111405} and Lemma \ref{le5.1}, $\mathbf{L}_t=(L_t^0,L_t,L_t^\ddagger)^\top$ satisfies
\begin{equation}\label{042803}
\begin{aligned}
d
\left(
\begin{array}{c}
{L}_t^0\\
L_t
\end{array}
\right)
=(\mathbb{D}_1+\mathbb{D}_2)^\top\left(
\begin{array}{c}
{L}_t^0\\
L_t
\end{array}
\right)dt
+ (\mathbb{F}_1+\mathbb{F}_2)^\top\left(
\begin{array}{c}
{L}_t^0\\
L_t
\end{array}
\right)dW_t^0,\
dL_t^\ddagger=f_2L_t^\ddagger dt+f_3L_t^\ddagger dW_t^0,
\end{aligned}
\end{equation}
with the initial condition $
\mathbf{L}_0=\rho(I+S_0\rho)^{-1}S_0\mathbf{L}_0+\rho (I+S_0\rho)^{-1}\Upsilon_0.
$
In particular, if $S_0=0$, then  $\mathbf{L}_0=\rho \Upsilon_0$; and if $\det(S_0^\top S_0)\neq 0$, then
$\mathbf{L}_0=(S_0^\top S_0)^{-1}S_0^\top(I+S_0\rho)S_0\rho(I+S_0\rho)^{-1}\Upsilon_0.$
Since (A6) can be readily verified, the following result follows directly from Theorem \ref{th5.1}.
\begin{corollary} Under (A5), assume $(S,\Upsilon,\upsilon)$ and $\mathbf{L}$ solve \eqref{042606}-\eqref{042607} and \eqref{042803}, respectively. Then an $\varepsilon_N$-Nash equilibrium $\Big(M_t^0,
M_t,\cdots,M_t\Big)$ for backward LQG-RMM is given by
\begin{equation}\label{050701}
\begin{aligned}
M_t^0=&{a_3^0}^\top \mathbf{L}_t
={a_3^0}^\top \Big(e^{\mathbb{H}t+\mathbb{K} W_t^0}(L_0^0,L_0),e^{(f_2-\frac12|f_3|^2)t+f_3W_t^0}
L_0^\ddagger\Big)^\top,\\
M_t=&(a_3+a_8)^\top \mathbf{L}_t=(a_3+a_8)^\top\Big(e^{\mathbb{H}t+\mathbb{K} W_t^0}(L_0^0,L_0),e^{(f_2-\frac12|f_3|^2)t+f_3W_t^0}
L_0^\ddagger\Big)^\top,
\end{aligned}
\end{equation}
with $\mathbb{H}=(\mathbb{D}_1+\mathbb{D}_2)^\top-\frac12(\mathbb{F}_1+\mathbb{F}_2)^\top
(\mathbb{F}_1+\mathbb{F}_2),\ \mathbb{K}=(\mathbb{F}_1+\mathbb{F}_2)^\top.$
\end{corollary}We have the following observations on \eqref{050701}:\\
(1) Both $M^{0}$ and $M$ are $\mathbb{F}^0$-adapted, and depend on the common noise $W^0$ only through $e^{\mathbb{K} W_t^0}$ and $e^{f_3W_t^0}$. If the BSDE's driver of \eqref{LQ-payoff} is independent to the intensity state $z$ (namely, $f_3^0=f_7^0=f_{3,7,11}=0$, so $\mathbb{K}=0$),
the equilibrium for each agent becomes deterministic.\\
(2) Unlike the forward LQG-RMM (Corollary \ref{le5.2}) with $\mathbb{F}^{i,0}$-adapted equilibrium, the equilibrium $(M_t)_{t\in[0,T]}$ in the backward case is $\mathbb{F}^{0}$-adapted, hence the idiosyncratic information driven by individual noise $\{W^i\}_{i=1}^N$ plays no role in the equilibrium. This is mainly because the driver of BSDE \eqref{LQ-payoff} is now linear and independent to the principal intensity state $Z^{i,i}$.\\
(3) If the major agent is absent, the equilibrium for each (minor) agent becomes $M_t=(a_3+a_8)^\top (0,0,L_t^\ddagger)^\top$. Comparing with $M_t$ in \eqref{050701}, one can see that the term $(a_3+a_8)^\top (L_t^0,L_t,0)^\top$ captures the influence of $\mathcal{A}_0$.

\begin{remark}\label{re5.1}
(i) When $f_3^0=f_3+f_{11},f_7^0=f_7=0$ ($\mathbb{K}^\top=f_3^0I$), \eqref{042606} can be solved as\begin{equation}\label{051403}
S_t=-\Big[(0,I_2)e^{
 \mathcal{B}(T-t)}(0,I)^\top\Big]^{-1}(0,I)e^{
 \mathcal{B}(T-t)}(I,0)^\top,\ t\in[0,T],
\end{equation}
where $\mathcal{B}=\left(
 \begin{array}{cc}
\widehat{\mathbb{C}}
&\widetilde{\mathbb{B}}+f_3^0(\mathbb{C}_5+\mathbb{C}_6)\rho\\
 -(\mathbb{C}_7+\mathbb{C}_8)&-\widetilde{\mathbb{C}}
 \end{array}
 \right)$, \ \ $\widehat{\mathbb{C}}=\mathbb{C}_3+\mathbb{C}_4+\rho(\mathbb{C}_7+\mathbb{C}_8)
+f_3^0(\mathbb{C}_5+\mathbb{C}_6)$
and (A5) is satisfied when $\det\{(0,I)e^{
 \mathcal{B}t}(0,I)^\top\}>0,\ \forall\ t\in[0,T].$

(ii) If $f_3^0=f_7^0=f_{3,7,11}=0$, $\mu^0_2=\mu_3=\mu_4=0$ and $f_8^0=f_8=f_{12}=0$, there has no weak-coupling by the control-average, and the $\varepsilon_N$-Nash equilibrium is given by
$M_t^0=-\frac12\big(R_0^{-1}f_4^0,R_0^{-1}f_8\big)
e^{\mathbb{H} t}(l_0^0,l_0)^\top,\
M_t=-\frac12 R^{-1}f_4e^{f_2t}
l_0^\ddagger.$ Besides, when the major agent is absent, the above result recovers that of \cite{HWW2016} (Theorem 3.1).
\end{remark}We present two concrete examples with more explicit representations for $\big(M_t^0,
M_t,\cdots,M_t\big)$.

\begin{example}\label{eg4}

Consider the backward LQG-RMM problem with: for $1\leq i\leq N,$
\begin{equation*}\label{eg1-payoff}
\begin{aligned}
J_0^N
=-\frac12|Y_0^0|^2-\frac12\mathbb{E}\Big[\int_0^T  |u_t^0-\mu^0_2\cdot u_t^{(N)}|^2dt\Big],\
J_{i}^N
=-\frac12|Y_0^i|^2
-\frac12\mathbb{E}\Big[\int_0^T |u_t^i-\mu_3\cdot u_t^{(N)}-\mu_4\cdot u_t^0|^2dt\Big],
\end{aligned}
\end{equation*}\begin{equation*}\label{eg1-BSDE}
\left\{
\begin{aligned}
-dY_t^0=&\Big[ Z_t^0+ \big(u_t^0
-\mu^0_2 u_t^{(N)}\big)\Big]dt-Z_t^0dW_t^0,\quad Y_T^0=\xi^0, \\
-dY_t^i=&\Big[ f_2Y_t^i
+{(1-f_2) Z_t^{i,0}}-f_2 Y_t^{(N)}+ f_{2}Z_t^{(N,0)}+\big(u_t^i
-\mu_4 u_t^{0}
-\mu_3 u_t^{(N)}\big)
\Big]dt\\
&-Z_t^{i,0}dW_t^0-Z_t^{i,i}dW_t^i,\quad
Y_T^i= \xi^i.
\end{aligned}
\right.
\end{equation*}Therefore, the major and each minor agent are weakly-coupled through their control-average along with an identical relative performance parameter. $\mathcal{B}$ in Remark \ref{re5.1} now reads as
$$\mathcal{B}=\left(
 \begin{array}{ccccc}
 0&\mu_4&-\mu^0_2&0&-\mu^0_2\\
 0&1&0&0&0\\
 0&0&0&0&0\\
 1&-\mu_4&\mu^0_2&1&\mu^0_2\\
 0&0&1&0&1\\
 \end{array}
 \right).$$
Since $\mathcal{B}^n=\mathcal{B}$, $n\geq 1$,
$\det\{(0,I)e^{\mathcal{B}t}(0,I)^\top\}=\det\left(
 \begin{array}{cc}
 e^t&\mu^0_2(e^t-1)\\
 0& e^t
 \end{array}
 \right)=e^{2t}>0, \forall\ t\in[0,T]$ hence (A5) holds. From \eqref{051403},
$$S_t=-\left(
 \begin{array}{ccc}
 1-e^{-(T-t)}&-\mu_4(1-e^{-(T-t)})&\mu^0_2e^{-(T-t)}(1-e^{-(T-t)})\\
 0&0&1-e^{-(T-t)}
 \end{array}
 \right),\ t\in[0,T].$$
We calculate from \eqref{050701} that, for $t\in[0,T],$
\begin{equation*}
\begin{aligned}
M_t^0=&-(1-\mu_3-\mu^0_2\mu_4)^{-1}e^{-\frac12 t+W_t^0}
\mathbb{E}\Big[e^{-\frac32 T+W_T^0}\Big((1-\mu_3)\xi^0+\mu^0_2(e^{f_2t}-1+e^{-T}+\mu_3-\mu_3e^{-T})\xi^1\Big)\Big],\\ M_t=&-(1-\mu_3-\mu^0_2\mu_4)^{-1}e^{-\frac12 t+W_t^0}
\mathbb{E}\Big[e^{-\frac32 T+W_T^0}\Big(\mu_4\xi^0+(e^{f_2t}-\mu_4\mu^0_2+\mu_4\mu^0_2e^{-T})\xi^1\Big)\Big].
\end{aligned}
\end{equation*}
\begin{remark}(1) The equilibrium $(M^0,M)$ have linear dependence on the terminal conditions $(\xi^0,\xi^1)$ via the parameters $(\mu_3,\mu^0_2,f_2,T,\mu_4)$.
(2) When $\mu^0_2=0$, the major's payoff is independent on the control-average of all minors, so her Nash strategy $M^0$ will not depend on $\xi^1$ any more, but may still get influence from each minor agent provided $\mu_3 \neq 0$.\end{remark}
\end{example}
\begin{example}\label{eg1}
Consider the backward LQG-RMM problem with:
\begin{equation*}\label{eg1-payoff}
\begin{aligned}
J_0^{N}
=-\gamma_0|Y_0^0|^2-\mathbb{E}\Big[\int_0^T  R_0|u_t^0|^2dt\Big],\
J_{i}^{N}
=-\gamma|Y_0^i|^2
-\mathbb{E}\Big[\int_0^T R|u_t^i|^2dt\Big],\ 1\leq i\leq N,
\end{aligned}
\end{equation*}
\begin{equation*}\label{eg1-BSDE}
\left\{
\begin{aligned}
-dY_t^0=&\Big(
f_6^0Y_t^{(N)}+f_8^0 u_t^{(N)}\Big)dt-Z_t^0dW_t^0,\quad Y_T^0=\xi^0,\\\
-dY_t^i=&\Big( f_2Y_t^i
+f_8 u_t^{0}
-f_{2} Y_t^{(N)}+f_{12} u_t^{(N)}
\Big)dt-Z_t^{i,0}dW_t^0-Z_t^{i,i}dW_t^i,\ Y_T^i= \xi^i,\ 1\leq i\leq N.
\end{aligned}
\right.
\end{equation*}
The major and each minor are interacted through the coupling $(Y^{(N)},u^{(N)})$. Moreover, each minor agent is also influenced by $u^0$ directly, the major's control. $\mathcal{B}$  in Remark \ref{re5.1} become
$$\mathcal{B}=\left(
 \begin{array}{ccccc}
 0&0&0&0&2\gamma_0f_6^0\\
 f_6^0&0&0&2\gamma_0f_6^0&0\\
 0&0&0&0&0\\
 0&0&0&0&-f_6^0\\
 0&0&0&0&0\\
 \end{array}
 \right).$$
 Now, $\mathcal{B}^2=0$,
and $\det\{(0,I_2)e^{\mathcal{B}t}(0,I_2)^\top\}=\det\left(
 \begin{array}{cc}
 1&-f_6^0t\\
 0&1
 \end{array}
 \right)=1>0$, $\forall\ t$, so (A5) holds.
Moreover, $S_t\equiv 0$ by \eqref{051403}, $M_t^0=-R_0^{-1}f_8f_6^0\gamma_0\Big(\mathbb{E}[\xi^0]+f_6^0T\mathbb{E}[\xi^1]\Big)t,\ M_t\equiv 0$ by \eqref{050701}.
\end{example}
\section{Conclusions} This paper studies a new class recursive major-minor (RMM) games featured by: (1) recursive functionals with nonlinear BSDE representations; (2) comprehensive and general weak-couplings. We propose a novel structural scheme to construct its auxiliary problem, a key step towards the desired $\varepsilon$-Nash equilibrium. In the RMM context, this scheme consists of a bilateral perturbation and a mixed triple-agent leader-follower-Nash analysis. In contrast to heuristic arguments in most MFG literature, such scheme indeed lay down an unified game-theoretic foundation to analyze more complex LP coupling structures, such as the ones with heterogenous robust beliefs, or with coalition interactions from nested information. We plan to address them in future.

\appendix
\section{Appendix}
\subsection{Proof of Theorem \ref{th4.1}}\label{app1}
For a fixed $N$, it suffices to verify the $\varepsilon_N$-Nash equilibrium property of $(u^{0,N},u^{1,N},\cdots$,
$u^{N,N})$ in side of $\mathcal{A}_0$. The verification in side of $\{\mathcal{A}_i\}_{i=1}^N$ is analogous thus we omit the details here. For this purpose, we consider the following limiting processes of \eqref{optimal-state} and \eqref{optimal-value}:
\begin{equation*}\label{limit-state}
\left\{
\begin{aligned}
dX_t^{0}=&b^0(t,X_t^{0},
u_t^{0,*};
 \overline{X}_t^{i},\overline{u}_t^{i,*})dt+{\sigma^{0}(t,{X}_t^{0}; \overline{X}_t^{i},\overline{u}_t^{i,*})}dW^0_t,\quad X_0^0=x_0^0,\\
dX_t^{i}=&
b(t,X_t^{i},{u}_t^{i,*};X_t^{0},u_t^{0,*};
\overline{X}_t^{i},\overline{u}_t^{i,*})dt
+\sigma(t,X_t^{i},X_t^{0};\overline{X}_t^{i},\overline{u}_t^{i,*})dW_t^i,\quad X_0^i=x_0,\\
-dY_t^{0}=&
f^0(t,\Theta_t^{0},
u_t^{0,*};  \overline{\Theta}_t^{i},\overline{u}_t^{i,*})dt-Z_t^0dW_t^0,\quad Y_T^{0}=\Phi^0(X_T^{0},\overline{X}_T^{i})+\xi^0, \\
-dY_t^{i}=&
f(t,\Theta_t^{i},{u}_t^{i,*};\Theta_t^{0},u_t^{0,*};\overline{\Theta}_t^{i},
\overline{u}_t^{i,*})dt
-Z_t^{i}d(W^0,W^i)_t,\quad Y_T^i=\Phi(X_T^{i},X_T^{0},\overline{X}_T^{i})+\xi^i,\\
\end{aligned}
\right.
\end{equation*}
and
$
J_0=\Gamma^0(Y_0^{0})+\mathbb{E}[
\int_0^Tg^0(t,\Theta_t^{0}, {u}_t^{0,*};\overline{\Theta}_t^{i},\overline{u}_t^{i,*})dt],
$
where $ u_t^{0,*}=\Psi^0(t,X_t^{0},
X_t^{1}),\ u_t^{i,*}=
\Psi(t,X_t^{0};X_t^{i}).$ By Burkholder-Davis-Gundy inequality and standard convergence estimates of SDEs (e.g., Theorem 10.1.7 in \cite{RR2006}), for $0\leq s\leq T,$ $j=0,1,\cdots,N,$ we have
\begin{equation*}\label{apr-1}
\begin{aligned}
&\mathbb{E}[\sup_{0\leq t\leq s}|X_t^{j,N}-X_t^j|^2]\leq
C\int_0^s\sum_{k=0,1,i}\mathbb{E}[\sup_{0\leq t\leq r}|X_r^{k,N}-X_r^k|^2]dr\\
&+C\int_0^s\mathbb{E}\Big[\Big|\frac{1}{N}\sum_{i=1}^NX_r^{i,N}-\mathbb{E}_r[X_r^{i}]\Big|^2\Big]dr
+C\int_0^s\mathbb{E}\Big[\Big|\frac{1}{N}\sum_{i=1}^N\Psi(r,X_r^{0,N},X_r^{i,N})-
\mathbb{E}_r[\Psi(r,X_r^0,X_r^{i})]
\Big|^2\Big]dr\\
&\leq C\int_0^s\mathbb{E}[\sup_{0\leq t\leq r}|X_r^{0,N}-X_r^0|^2]dr
+C\int_0^s\sup_{1\leq i\leq N}\mathbb{E}[\sup_{0\leq t\leq r}|X_r^{i,N}-X_r^i|^2]dr+\frac{C}{N}.
\end{aligned}
\end{equation*}
It follows from Gronwall inequality that
\begin{equation}\label{apr-3}
\sup_{0\leq i\leq N}\mathbb{E}[\sup_{0\leq t\leq T}|X_t^{i,N}-X_t^i|^2]\leq  {\frac{C}{N}}.
\end{equation}
Applying It\^o's formula to $|Y_t^{0,N}-Y_t^0|^2$ and $|Y_t^{i,N}-Y_t^i|^2$, we have
\begin{small}
\begin{equation}\nonumber
\begin{aligned}
&\mathbb{E}\big[|Y_t^{j,N}-Y_t^j|^2+\int_t^T|Z_s^{j,N}-Z_s^j|^2ds\big]\leq
C\int_t^T\mathbb{E}[|Y_s^{j,N}-Y_s^j|^2]ds+C\mathbb{E}[ |X_T^{j,N}-X_T^j|^2]
+C_j\mathbb{E}[ |X_T^{0,N}-X_T^0|^2]\\
&+C\mathbb{E}\big[\big|\frac{1}{N}\sum_{i=1}^NX_T^{i,N}
-\mathbb{E}_T[X_T^{i}]
\big|^2\big] +C\int_t^T\sum_{k=0,1,i}\mathbb{E}[|X_s^{k,N}-X_s^k|^2]ds
+C_j\mathbb{E}\big[\int_t^T(|Y_s^{0,N}-Y_s^0|^2+|Z_s^{0,N}-Z_s^0|^2)ds\big]\\
\end{aligned}
\end{equation}
\end{small}
\begin{small}
\begin{equation}\label{apr-5}
\begin{aligned}
&+C\int_t^T\mathbb{E}\Big[\Big|\frac{1}{N}\sum_{i=1}^NX_s^{i,N}
-\mathbb{E}_s[X_s^{i}]\Big|^2\Big]ds
+C\int_t^T\mathbb{E}\Big[\Big|\frac{1}{N}\sum_{i=1}^N\Psi(r,X_r^{0,N},X_r^{i,N})-
\mathbb{E}_r[\Psi(r,X_r^0,X_r^{i})]
\Big|^2\Big]dr\\
&+C\int_t^T\mathbb{E}\Big[\Big|\frac{1}{N}\sum_{i=1}^NY_s^{i,N}
-\mathbb{E}_s[Y_s^{i}]\Big|^2\Big]ds+C\int_t^T\mathbb{E}\Big[\Big|\frac{1}{N}\sum_{i=1}^NZ_s^{i,N}
-\mathbb{E}_s[Z_s^{i}]\Big|^2\Big]ds\\
&\leq C\int_t^T\mathbb{E}[|Y_s^{i,N}-Y_s^i|^2]ds+C\mathbb{E}[ |X_T^{j,N}-X_T^j|^2] +C_j\mathbb{E}[ |X_T^{0,N}-X_T^0|^2]+C\int_t^T\sum_{k=0,1,i}\mathbb{E}[|X_s^{k,N}-X_s^k|^2]ds\\
&+C_j\mathbb{E}\big[\int_t^T(|Y_s^{0,N}-Y_s^0|^2+|Z_s^{0,N}-Z_s^0|^2)ds\big]
+\frac{C}{N},
\end{aligned}
\end{equation}
\end{small}
where $C_0=0$ and $\{C_j\}_{j=1}^N$ stands for some positive constant.
Combining \eqref{apr-3} and \eqref{apr-5},
\begin{equation*}\label{apr-6}
\sup_{0\leq i\leq N}\Big\{\sup_{0\leq t\leq T}\mathbb{E}[|Y_t^{i,N}-Y_t^i|^2]+\mathbb{E}[\int_0^T|Z_t^{i,N}-Z_t^i|^2dt]\Big\}
\leq \frac{C}{N}.
\end{equation*}As to the functionals,\begin{small}
\begin{equation}\nonumber
\begin{aligned}
&|J^N_0-J_0|\leq |\Gamma^0(Y_0^{0,N})-\Gamma^0(Y_0^0)|
+
\mathbb{E}\Big[
\int_0^T\Big|g^0\Big(t,\Theta_t^{0,N}, \Psi^0(t,X_t^{0,N},
X_t^{1,N});\Theta_t^{(N)},u_t^{N}\Big)\\
&-g^0\Big(t,\Theta_t^{0},  \Psi^0(t,X_t^{0},
X_t^{1});\overline{\Theta}_t^{i},\mathbb{E}_t\big[ \Psi(t,X_t^{0},
X_t^{i})\big]\Big)\Big|dt\Big]\\
&\leq C(1+|Y_0^{0,N}|+|Y_0^0|)\cdot |Y_0^{0,N}-Y_0^{0}|\\
&+C\mathbb{E}\Big[
\int_0^T\big(1+|\Theta_t^{0,N}|+|\Theta_t^{0}|+|\Theta_t^{(N)}|+|\overline{\Theta}_t^{i}|+
|u_t^{(N)}|+|\mathbb{E}_t[ \Psi(t,X_t^0,X_t^1)]|\big)\\
&\big(|\Theta_t^{0,N}-\Theta_t^{0}|
+|\Theta_t^{(N)}-\overline{\Theta}_t^{i}|+|u_t^{(N)}-\mathbb{E}_t[ \Psi(t,X_t^0,X_t^1)]|
\big)dt\\
&\leq C(1+|Y_0^{0,N}|^2+|Y_0^0|^2)\cdot |Y_0^{0,N}-Y_0^{0}|^2\\
&+C\Big(\mathbb{E}\Big[
\int_0^T\big(1+|\Theta_t^{0,N}|^2+|\Theta_t^{0}|^2+|\Theta_t^{(N)}|^2+|\overline{\Theta}_t^{i}|^2+
|u_t^{(N)}|^2+| \Psi(t,X_t^0,X_t^1)|^2\big)dt\Big)^{\frac12}\cdot\\
&\Big(\mathbb{E}\Big[
\int_0^T\big(|\Theta_t^{0,N}-\Theta_t^{0}|^2
+|\Theta_t^{(N)}-\overline{\Theta}_t^{i}|^2+
|u_t^{(N)}-\mathbb{E}_t[ \Psi(t,X_t^0,X_t^1)]|^2
\big)dt\Big)^{\frac12}
\leq \frac{C}{\sqrt{N}}.
\end{aligned}
\end{equation}
\end{small}We study the uniliteral deviation of $\mathcal{A}_0$ from the strategy $u_t^{0,N}=\Psi^0(t,X_t^{0,N},
X_t^{1,N})$. Assume now that $\mathcal{A}_0$ adopts a different control $u^0\in\mathcal{U}_d^0$
and $\{\mathcal{A}_i\}_{i\geq 1}^N$ keep to apply $\{u_t^{i,N}\}_{i\geq 1}^N$. The resulting perturbed states, denoted by $(\widehat{X}_t^{i,N})_{i\geq 0}$, should satisfy
\begin{equation*}\label{apr-8}
\left\{
\begin{aligned}
d\widehat{X}_t^{0,N}=&b^0\big(t,\widehat{X}_t^{0,N},
u^0_t; \widehat{X}_t^{(N)},{u_t^{(N)}}\big)dt+{\sigma^{0}\big(t,\widehat{X}_t^{0,N}; \widehat{X}_t^{(N)},{u_t^{(N)}}\big)}dW^0_t,\quad
\widehat{X}_0^{0,N}=x_0^0,\\
d\widehat{X}_t^{i,N}=&b\big(t,\widehat{X}_t^{i,N},
 u_t^{i,*,N};
\widehat{X}_t^{0,N},u_t^0;\widehat{X}_t^{(N)},u_t^{(N)}\big)dt
+\sigma\big(t,\widehat{X}_t^{i,N},\widehat{X}_t^{0,N};
\widehat{X}_t^{(N)},u_t^{(N)}\big)dW_t^i,\quad
\widehat{X}_0^{i,N}=x_0,\\
\end{aligned}
\right.
\end{equation*}
and the related limiting processes is given by
\begin{equation*}\label{apr-9}
\left\{
\begin{aligned}
d\widehat{X}_t^{0}=&b^0\big(t,\widehat{X}_t^{0},
u^0_t; \mathbb{E}_t[\widehat{X}_t^{i}],\overline{u}_t^{i,N}\big)dt
+{\sigma^{0}\big(t,\widehat{X}_t^{0}; \mathbb{E}_t[\widehat{X}_t^{i}],\overline{u}_t^{i,N}\big)}dW^0_t,\quad
\widehat{X}_0^0=x_0^0,\\
d\widehat{X}_t^{i}=&b\big(t,\widehat{X}_t^{i},
{u}_t^{i,N};
\widehat{X}_t^{0},u_t^0;\mathbb{E}_t[\widehat{X}_t^{i}],
\overline{u}_t^{i,N}\big)dt
+\sigma\big(t,\widehat{X}_t^{i},\widehat{X}_t^{0};\mathbb{E}_t[\widehat{X}_t^{i}],
\overline{u}_t^{i,N}\big)dW_t^i,\quad
\widehat{X}_0^i=x_0.
\end{aligned}
\right.
\end{equation*}
Similar to \eqref{apr-3},
$
\sup_{0\leq i\leq N}\mathbb{E}[\sup_{0\leq t\leq T}|\widehat{X}_t^{i,N}-\widehat{X}_t^i|^2]\leq \frac{C}{N}.
$
By the same estimates as in \eqref{apr-7},
\begin{equation}\label{apr-11}
|\widehat{J}^N_0-\widehat{J}_0|\leq \frac{C}{\sqrt{N}}.
\end{equation}
Since $\big(
 u^{0,*}_t,
 u^{1,*}_t, u_t^{j,*}\big)$ is an equilibrium  strategy of the limiting triple-agent game problem, it is clear that $\widehat{J}_0\leq {J}_0$, and combining \eqref{apr-7} and \eqref{apr-11}, we get the desired result for $\mathcal{A}_0$.

\subsection{Some notations in Section \ref{LQ}}\label{app2}
1. Constants in \eqref{11150101}: $a=\frac{1}{2}(1-\mu_3-\mu^0_2\mu_4)^{-1}$ and
\begin{small}
\begin{equation}\nonumber
\begin{array}{llll}
& {a_1^0}^\top=a\big((1-\mu_3)R_0^{-1}b_2^0,(1-\mu_3)R_0^{-1}b_4,\mu^0_2 R^{-1}b_2\big),
& {a_3^0}^\top=-a\big((1-\mu_3)R_0^{-1}f_4^0,(1-\mu_3)R_0^{-1}f_8,\mu^0_2 R^{-1}f_4\big),
\\
& a_1^\top=a\big(\mu_4R_0^{-1}b_2^0,\mu_4R_0^{-1}b_4,(\mu_3+\mu_4\mu^0_2) R^{-1}b_2\big),
& a_3^\top=-\big(a\mu_4R_0^{-1}f_4^0,a\mu_4R_0^{-1}f_8,(a-\frac12) R^{-1}f_4\big),
\\
& a_7^\top=\frac12 R^{-1}b_2(0,0,1),
& a_8^\top=-\frac12 R^{-1}f_4(0,0,1).
\\
\end{array}
\end{equation}
\end{small}
2. Matrices in equation \eqref{043001}:
\begin{small}
$$\mathbb{A}_1=
\left(
\begin{array}{cc}
b_1^0&0 \\
b_3&b_1\\
\end{array}
\right),\
\mathbb{A}_2=
\left(
\begin{array}{cc}
0&b_3^0 \\
0&b_5\\
\end{array}
\right),\
\mathbb{A}_3=
\left(
\begin{array}{cc}
f_1^0&0\\
f_5&f_1\\
\end{array}
\right),\
\mathbb{A}_4=
\left(
\begin{array}{cc}
0&f_5^0\\
0&f_9\\
\end{array}
\right),\
\mathbb{A}_5=
\left(
\begin{array}{cc}
-2Q_0&0\\
2\mu^0_1Q_0&0\\
2\mu_2Q&-2Q\\
\end{array}
\right),\
$$
$$\mathbb{A}_6=
\left(
\begin{array}{cc}
0&2\mu^0_1Q_0\\
0&-2|\mu^0_1|^2Q_0\\
0&2\mu_1Q\\
\end{array}
\right),\
\mathbb{A}_7=
\left(
\begin{array}{cc}
\Phi_1^0&0\\
\Phi_2&\Phi_1\\
\end{array}
\right),\
\mathbb{A}_8=
\left(
\begin{array}{ccccc}
0&\Phi_2^0\\
0&\Phi_3\\
\end{array}
\right),\
\mathbb{B}_1=
\left(
\begin{array}{c}
0\\
b_2a_7^\top\\
\end{array}
\right),\
$$
$$
\mathbb{B}_2=
\left(
\begin{array}{c}
b_2^0{a_1^0}^\top+b_4^0{a}_1^\top+b_4^0{a}_7^\top\\
(b_2+b_6)a_1^\top+b_4{a_1^0}^\top+b_6{a}_7^\top\\
\end{array}
\right),\
\mathbb{B}_3=
\left(
\begin{array}{c}
0\\
f_4a_7^\top\\
\end{array}
\right), \
\mathbb{B}_4=
\left(
\begin{array}{c}
f_4^0{a_1^0}^\top+f_8^0(a_1+a_7)^\top\\
f_4a_1^\top+f_8{a_1^0}^\top+f_{12}(a_1+a_7)^\top\\
\end{array}
\right), \
$$
$$
\mathbb{B}_5=
\left(
\begin{array}{ccc}
b_1^0&b_3&0\\
0&b_1&0\\
0&0&b_1
\end{array}
\right), \
\mathbb{B}_6=
\left(
\begin{array}{ccc}
0&0&0\\
b_3^0&b_5&0\\
0&0&0
\end{array}
\right), \
\mathbb{C}_1=
\left(
\begin{array}{c}
0\\
b_2a_8^\top\\
\end{array}
\right),\
\mathbb{C}_2=
\left(
\begin{array}{c}
b_2^0{a_3^0}^\top+b_4^0{a}_3^\top+b_4^0{a}_8^\top\\
(b_2+b_6)a_3^\top+b_4{a_3^0}^\top+b_6{a}_8^\top\\
\end{array}
\right),\
$$
$$
\mathbb{C}_3=
\left(
\begin{array}{ccc}
f_2^0&f_6&0\\
0&f_2&0\\
0&0&f_2\\
\end{array}
\right), \
\mathbb{C}_4=
\left(
\begin{array}{ccccc}
0&0&0\\
f_6^0&f_{10}&0\\
0&0&0\\
\end{array}
\right), \
\mathbb{C}_5=
\left(
\begin{array}{ccc}
f_3^0&f_7&0\\
0&f_3&0\\
0&0&f_3\\
\end{array}
\right), \
\mathbb{C}_6=
\left(
\begin{array}{ccccc}
0&0&0\\
f_7^0&f_{11}&0\\
0&0&0\\
\end{array}
\right), \
$$
$$
\mathbb{C}_7=
\left(
\begin{array}{c}
0\\
f_4a_8^\top\\
\end{array}
\right),\
\mathbb{C}_8=
\left(
\begin{array}{c}
f_4^0{a_3^0}^\top+f_8^0(a_3+a_8)^\top\\
f_4a_3^\top+f_8{a_3^0}^\top+f_{12}(a_3+a_8)^\top\\
\end{array}
\right),\
\mathbb{C}_9=
\left(
\begin{array}{ccccc}
-f_1^0&-f_5&0\\
0&-f_{1}&0\\
0&0&-f_1\\
\end{array}
\right), \
$$
$$
\mathbb{C}_{10}=
\left(
\begin{array}{ccccc}
0&0&0\\
-f_5^0&-f_{9}&0\\
0&0&0\\
\end{array}
\right), \
\mathbb{C}_{11}=
\left(
\begin{array}{ccc}
-\Phi_1^0&-\Phi_2&0\\
0&-\Phi_1&0\\
0&0&-\Phi_1\\
\end{array}
\right),\
\mathbb{C}_{12}=
\left(
\begin{array}{ccc}
0&0&0\\
-\Phi_2^0&-\Phi_3&0\\
0&0&0\\
\end{array}
\right),\
$$
$$
\mathbb{D}_1=
\left(
\begin{array}{cc}
f_2^0&0\\
f_6&f_2\\
\end{array}
\right),\
\mathbb{D}_2=
\left(
\begin{array}{cc}
0&f_6^0\\
0&f_{10}\\
\end{array}
\right),\
\mathbb{F}_1=
\left(
\begin{array}{cc}
f_3^0&0\\
f_7&f_3\\
\end{array}
\right),\
\mathbb{F}_2=
\left(
\begin{array}{cc}
0&f_7^0 \\
0&f_{11}\\
\end{array}
\right),\
 \mathbb{\rho}=
\left(
\begin{array}{cc}
2\gamma_0&0\\
0&0\\
0& 2\gamma\\
\end{array}
\right).
 $$
\end{small}
3. The $5\times 5$ block matrices in equations \eqref{31204} and \eqref{31205}:
$$\mathbb{H}_1=
\left(
\begin{array}{cc}
\mathbb{A}_1+\mathbb{A}_2& \mathbb{C}_1+\mathbb{C}_2\\
\rho(\mathbb{A}_3+\mathbb{A}_4)&\rho(\mathbb{C}_7+\mathbb{C}_8)\\
\end{array}
\right),\
\mathbb{H}_2=
\left(
\begin{array}{cc}
(\mathbb{C}_7+\mathbb{C}_8)\rho+\mathbb{D}_1+\mathbb{D}_2
&\mathbb{B}_3+\mathbb{B}_4 \\
(\mathbb{C}_9+\mathbb{C}_{10})\rho&\mathbb{B}_5+\mathbb{B}_6\\
\end{array}
\right),\
$$
$$
\mathbb{H}_3=
\left(
\begin{array}{cc}
(\mathbb{C}_1+\mathbb{C}_2)\rho
&\mathbb{B}_1+\mathbb{B}_2 \\
(\mathbb{C}_3+\mathbb{C}_{4})\rho+\rho(\mathbb{C}_7+\mathbb{C}_{8})\rho
+\rho(\mathbb{D}_1+\mathbb{D}_{2})&\rho(\mathbb{B}_3+\mathbb{B}_4)\\
\end{array}
\right),\
\mathbb{H}_4=
\left(
\begin{array}{cc}
\mathbb{F}_1+\mathbb{F}_2
&0 \\
0&0\\
\end{array}
\right),\
$$
$$
\mathbb{H}_5=
\left(
\begin{array}{cc}
0&0 \\
0&\mathbb{C}_5+\mathbb{C}_6\\
\end{array}
\right),\
\mathbb{H}_6=
\left(
\begin{array}{cc}
0&0 \\
(\mathbb{C}_5+\mathbb{C}_6)\rho&0\\
\end{array}
\right),\
\mathbb{H}_7=
\left(
\begin{array}{cc}
\mathbb{A}_3+\mathbb{A}_4&\mathbb{C}_7+\mathbb{C}_8 \\
\mathbb{A}_5+\mathbb{A}_6&\mathbb{C}_9+\mathbb{C}_{10}\\
\end{array}
\right),\
$$
$$
\mathbb{G}_1=
\left(
\begin{array}{cc}
0&0 \\
(\mathbb{C}_{11}+\mathbb{C}_{12})\rho&0\\
\end{array}
\right),\
\mathbb{G}_2=
\left(
\begin{array}{cc}
\mathbb{A}_7+\mathbb{A}_8&0 \\
0&\mathbb{C}_{11}+\mathbb{C}_{12}\\
\end{array}
\right),\
\widetilde{\rho}=
\left(
\begin{array}{cc}
0&0 \\
\rho&0\\
\end{array}
\right).\
$$
4. Constant vectors in equation \eqref{042601}
\begin{small}
$$
\begin{aligned}
\mathbf{\Lambda}_1=&\big(b_3,b_5,
-[(b_2+b_6)\mu_4+b_4(1-\mu_3)]aR_0^{-1}f_4^0,
-[(b_2+b_6)\mu_4+b_4(1-\mu_3)]aR_0^{-1}f_8,
-[b_2+b_4\mu^0_2+b_6]aR^{-1}f_4\big),\\
\mathbf{\Lambda}_2=&\big(0,0,
\mathbf{\Lambda}_3\big),\
\mathbf{\Lambda}_4=(2Q \mu_2,2Q\mu_1,0,0,-f_1),\\
\mathbf{\Lambda}_3=&\big(
[b_2\mu_4+b_4(1-\mu_3)+b_6\mu_4]aR_0^{-1}b_2^0,
[b_2\mu_4+b_4(1-\mu_3)+b_6\mu_4]aR_0^{-1}b_4,
[b_2(\mu_3+\mu_4\mu^0_2)+b_4\mu^0_2+b_6]aR^{-1}b_2\big).
\end{aligned}
$$
\end{small}

\subsection{Proof of Proposition \ref{LQ-PRO}}
The uniqueness part follows by the standard arguments. Thus, we focus only on the existence part. By Lemma 6.1, we can substitute the solution $(\overline{\mathbf{X}}_t,\overline{\mathbf{L}}_t)^\top=(\overline{X}^0_t,\overline{X}^1_t,\overline{L}_t^0,
\overline{L}_t,\overline{L}^\ddagger_t)^\top,$ $(\overline{\mathbf{Y}}_t,\overline{\mathbf{P}}_t)^\top=(\overline{Y}_t^0,\overline{Y}_t^1,\overline{P}^0_t,
\overline{P}_t,\overline{P}^\ddagger_t)^\top,\ (\overline{\mathbf{Z}}_t,\overline{\mathbf{Q}}_t)^\top=(\overline{Z}_t^0,\overline{Z}_t^1,\overline{Q}^{00}_t,
\overline{Q}^{10}_t,\overline{Q}^{10,\ddagger}_t)^\top$ of \eqref{111405} into \eqref{043001}. Then, its solution can be constructed by the following three steps.

\emph{{Step 1.} The construction of $(X^0,Y^0,Z^0)$ and $(L^0,L,L^\ddagger;P^0,Q^{00}, Q^{01};P,Q^{10},Q^{11})$.}

The triple $(X^0,Y^0,Z^0)$ (if exists), should be $\mathbb{F}^0$-adapted, and define
$X^0_t=\overline{X}^0_t,\ Y^0_t=\overline{Y}^0_t,\ Z^0_t=\overline{Z}^0_t.$
From \eqref{043001}, the adjoint equations (except those of $(P^\ddagger,Q^{10,\ddagger},Q^{11,\ddagger}$)) satisfy
\begin{equation*}\label{pr1}
\left\{
\begin{aligned}
dL_t^0=&[f_2^0 L_t^0+f_6 L_t]dt+[f_3^0 L_t^0+f_7 L_t]dW^0_t,\quad L_0^0=2\gamma_0Y_0^0,\\
dL_t=&\big[f_2 L_t+f_{10} \overline{L}_t
+f_6^0 \overline{L}_t^0\big]dt
+\big[f_3 L_t+f_{11} \overline{L}_t
+f_7^0 \overline{L}_t^0\big]dW_t^0,\quad  L_0=0,\\
dL_t^{\ddagger}=&f_2 L_t^\ddagger dt
+f_3 L_t^\ddagger dW_t^0,\quad \ L_0^\ddagger=2\gamma \overline{Y}_0^1,\\
\end{aligned}
\right.
\end{equation*}
\begin{equation}\label{pr1}
\left\{
\begin{aligned}
-d{P_t^0}=&\Big[b_1^0{P_t^0}
+b_3{P_t}-f_1^0 L_t^0-f_5 L_t
-2Q_0\big(
{X_t^0}-\mu^0_1\overline{X}_t^1\big)\Big]dt
-{Q_t^{00}}dW^0_t
-{Q_t^{01}}dW_t^1,\\
-d{P_t}=&
\Big[b_1{P_t}-f_1 L_t
+b_5\overline{P}_t
+b_3^0\overline{P}_t^0-f_5^0\overline{L}_t^0
-f_9\overline{L}_t
+2\mu^0_1Q_0\big(X_t^0-\mu^0_1\overline{X}_t^1\big)\Big]dt\\
&-{Q_t^{10}}dW^0_t-{Q_t^{11}}dW_t^1,\\
{P_T^0}=&-L_T^0\Phi_1^0-L_T\Phi_2,\quad
{P_T}=-L_T\Phi_1-\overline{L}_T^0\Phi_2^0-\overline{L}_T\Phi_3.
\end{aligned}
\right.
\end{equation}
\eqref{pr1} is a linear decoupled FBSDE with $L^{2}$ nonhomogeneous terms, thus it admits a unique solution $(L^0,L,L^\ddagger;P^0,Q^{00}, Q^{01};P,Q^{10},Q^{11})$. Moreover, comparing \eqref{pr1}, \eqref{111405}, we have$$L_t^0=\overline{L}_t^0,\ L_t=\overline{L}_t,\ L_t^{\ddagger}=\overline{L}^\ddagger_t,\ P_t^0=\overline{P}_t^0,\ Q^{00}_t=
\overline{Q}_t^{00}, Q^{01}_t=0,\ P_t=\overline{P}_t,\ Q^{10}_t=
\overline{Q}_t^{10}, Q^{11}_t=0,\ t\in[0,T].$$

\emph{Step 2. The construction of $(X^1,P^\ddagger,Q^{10,\ddagger},Q^{11,\ddagger})$.}

It follows from the system \eqref{043001} that the 4-tuple $(X^1,P^\ddagger,Q^{10,\ddagger},Q^{11,\ddagger})$ satisfies
\begin{equation}\label{pr2}
\left\{
\begin{aligned}
dX_t^1=&\Big[b_1X_t^1+b_2 a_7 P_t^\ddagger
+\tilde{b}_t\Big]dt
+\sigma dW_t^1,\ {{X}_0^{1}}=x_0,\\
-d{{P}_t^\ddagger}=&\Big[b_1{P_t^\ddagger}-2Q{X_t^1}
+\tilde{f}_t\Big]dt-{Q_t^{10,\ddagger}}dW^0_t
-{Q_t^{11,\ddagger}}dW_t^1,\
{P_T^\ddagger}=
-L_T^\ddagger\Phi_1,
\end{aligned}
\right.
\end{equation}
where $
\tilde{b}_t=\mathbf{\Lambda}_1\overline{\mathbf{X}}_t
+\mathbf{\Lambda}_2\overline{\mathbf{Y}}_t
=[\mathbf{\Lambda}_1+\mathbf{\Lambda}_2(I+S_t\rho)^{-1}S_t]\overline{\mathbf{X}}_t
+\mathbf{\Lambda}_2(I+S_t\rho)^{-1}\Upsilon_t$ and $
\tilde{f}_t=\mathbf{\Lambda}_4\overline{\mathbf{X}}_t$
 are $\mathbb{F}^0$-adapted.
 Next
we introduce the following Riccati equation and BSDE
\begin{equation}\nonumber
\left\{
\begin{aligned}
&\dot{\Sigma}_t+2b_1\Sigma_t+\frac12 R^{-1}|b_2|^2|\Sigma_t|^2-2Q=0,\ \ \  \Sigma_T=0,\\
&-dp_t=[(\frac12 R^{-1}|b_2|^2\Sigma_t+b_1)p_t+\Sigma_t\tilde{b}_t+\tilde{f}_t]dt-q_t^0dW_t^0-q_t^1dW_t^1,\ \ \  p_T=-L_T^\ddagger\Phi_1,
\end{aligned}
\right.
\end{equation}
Under (A6), it follows from Theorem 4.3 in \cite{MY1999} (see {p.48}) that $\Sigma$ takes the form \eqref{042503}.
Notice that $\Sigma$ is bounded, the above linear BSDE has an unique $\mathbb{F}^0$-adapted solution $(p,q^0,0)$ with \eqref{042601}.
By the relation $
P_t^\ddagger=\Sigma_t X_t^{1}+p_t,$ one can show the well-posedness of \eqref{pr2}.

\emph{Step 3. The construction of $(Y^1,Z^{1,0},Z^{1,1})$.}

Denote $
f_t=f_1X_t^1+f_5 X_t^{0}+f_6 Y_t^{0}+f_7 Z_t^{0}
+f_9 \overline{X}_t^{1}+f_{10} \overline{Y}_t^{1}+{f_{11} \overline{Z}_t^{1,0}}+f_4 a_7 P_t^\ddagger
+(f_4 a_1+f_8a_1^0+f_{12}a_1+f_{12}a_7)^\top
 \overline{\mathbf{P}}_t
+[(f_4+f_{12})( a_3+a_8)+f_{8}{a}_3^0]^\top \mathbf{L}_t
+f_4 a_7^\top\mathbf{P}_t.
$
Then $(Y^1,Z^{1,0},Z^{1,1})$ is the unique $\mathbb{F}^{0,1}$-adapted solution of the following BSDE
\begin{equation*}\label{pr3}
\begin{aligned}
-dY_t^1=\Big[f_2Y_t^1+{f_3Z_t^{1,0}}+f_t
\Big]dt-{Z_t^{1,0}dW_t^0-Z_t^{1,1}dW_t^1},\
Y_T^1= \Phi_1X_T^1 +\Phi_2X_T^0 +\Phi_3\overline{X}_T^{1}+\xi^1.
\end{aligned}
\end{equation*}Finally, combining the above three steps we construct a solution of
the system \eqref{043001}.

\bibliographystyle{SIAM}
\small\vspace{-2 mm}
\bibliography{mybib}

\begin{thebibliography}{10}

\bibitem{AB1984}
{\sc A.~Agrawal and R.~E. Barlow}, {\em A survey of network reliability and
  domination theory}, Operations Research, 32 (1984), pp.~478--492.

\bibitem{AD2011}
{\sc D.~Andersson and B.~Djehiche}, {\em A maximum principle for {SDEs} of
  mean-field type}, Applied Mathematics \& Optimization, 63 (2011),
  pp.~341--356.

\bibitem{AMS1995}
{\sc R.~J. Aumann, M.~Maschler, and R.~E. Stearns}, {\em Repeated games with
  incomplete information}, MIT press, 1995.

\bibitem{AH2021}
{\sc J.~Aurand and Y.-J. Huang}, {\em Mortality and healthcare: A stochastic
  control analysis under {Epstein--Zin} preferences}, SIAM Journal on Control
  and Optimization, 59 (2021), pp.~4051--4080.

\bibitem{BCY2016}
{\sc A.~Bensoussan, M.~H. Chau, and S.~C. Yam}, {\em Mean field games with a
  dominating player}, Applied Mathematics \& Optimization, 74 (2016),
  pp.~91--128.

\bibitem{BCR2024}
{\sc P.~Bergault, P.~Cardaliaguet, and C.~Rainer}, {\em Mean field games in a
  stackelberg problem with an informed major player}, SIAM Journal on Control
  and Optimization, 62 (2024), pp.~1737--1765.

\bibitem{BLP2014}
{\sc R.~Buckdahn, J.~Li, and S.~Peng}, {\em Nonlinear stochastic differential
  games involving a major player and a large number of collectively acting
  minor agents}, SIAM Journal on Control and Optimization, 52 (2014),
  pp.~451--492.

\bibitem{CCP2020}
{\sc P.~Cardaliaguet, M.~Cirant, and A.~Porretta}, {\em Remarks on {Nash}
  equilibria in mean field game models with a major player}, Proceedings of the
  American Mathematical Society, 148 (2020), pp.~4241--4255.

\bibitem{CZ2018}
{\sc R.~Carmona, F.~Delarue, R.~Carmona, and F.~Delarue}, {\em Extensions for
  volume {I}}, Probabilistic Theory of Mean Field Games with Applications I:
  Mean Field FBSDEs, Control, and Games,  (2018), pp.~619--680.

\bibitem{CW2017}
{\sc R.~Carmona and P.~Wang}, {\em An alternative approach to mean field game
  with major and minor players, and applications to herders impacts}, Applied
  Mathematics \& Optimization, 76 (2017), pp.~5--27.

\bibitem{CZ2016}
{\sc R.~A. Carmona and X.~Zhu}, {\em A probabilistic approach to mean field
  games with major and minor players}, Annals of Applied Probability, 26
  (2016), pp.~1535--1580.

\bibitem{CE2002}
{\sc Z.~Chen and L.~Epstein}, {\em Ambiguity, risk, and asset returns in
  continuous time}, Econometrica, 70 (2002), pp.~1403--1443.

\bibitem{D2002}
{\sc F.~Delarue}, {\em On the existence and uniqueness of solutions to {FBSDEs}
  in a non-degenerate case}, Stochastic Processes and Their Applications, 99
  (2002), pp.~209--286.

\bibitem{DHW2018}
{\sc K.~Du, J.~Huang, and Z.~Wu}, {\em Linear quadratic mean-field-game of
  backward stochastic differential systems}, Mathematical Control and Related
  Fields, 8 (2018), pp.~653--678.

\bibitem{DE1992}
{\sc D.~Duffie and L.~G. Epstein}, {\em Stochastic differential utility},
  Econometrica, 60 (1992), pp.~353--394.

\bibitem{EPQ2001}
{\sc N.~El~Karoui, S.~Peng, and M.~C. Quenez}, {\em A dynamic maximum principle
  for the optimization of recursive utilities under constraints}, Annals of
  Applied Probability, 11 (2001), pp.~664--693.

\bibitem{FHH2022}
{\sc X.~Feng, Y.~Hu, and J.~Huang}, {\em Backward stackelberg differential game
  with constraints: a mixed terminal-perturbation and linear-quadratic
  approach}, SIAM Journal on Control and Optimization, 60 (2022),
  pp.~1488--1518.

\bibitem{HL2019}
{\sc Z.~Hellman and Y.~J. Levy}, {\em Measurable selection for purely atomic
  games}, Econometrica, 87 (2019), pp.~593--629.

\bibitem{HJX2023}
{\sc M.~Hu, S.~Ji, and X.~Xue}, {\em Optimization under rational expectations:
  A framework of fully coupled forward-backward stochastic linear quadratic
  systems}, Mathematics of Operations Research, 48 (2023), pp.~1767--1790.

\bibitem{HHL2022}
{\sc Y.~Hu, J.~Huang, and W.~Li}, {\em Backward stochastic differential
  equations with conditional reflection and related recursive optimal control
  problems}, SIAM Journal on Control and Optimization, 62 (2024),
  pp.~2557--2589.

\bibitem{HLZ2023}
{\sc J.~Huang, W.~Li, and H.~Zhao}, {\em A class of optimal control problems of
  forward--backward systems with input constraint}, Journal of Optimization
  Theory and Applications, 199 (2023), pp.~1050--1084.

\bibitem{HWW2016}
{\sc J.~Huang, S.~Wang, and Z.~Wu}, {\em Backward mean-field
  linear-quadratic-{Gaussian (LQG)} games: full and partial information}, IEEE
  Transactions on Automatic Control, 61 (2016), pp.~3784--3796.

\bibitem{H2010}
{\sc M.~Huang}, {\em Large-population {LQG} games involving a major player: the
  {N}ash certainty equivalence principle}, SIAM Journal on Control and
  Optimization, 48 (2010), pp.~3318--3353.

\bibitem{HMC2006}
{\sc M.~Huang, R.~P. Malham\'e, and P.~E. Caines}, {\em Large population
  stochastic dynamic games: closed-loop mckean-vlasov systems and the nash
  certainty equivalence principle}, Communications in Information and Systems,
  6 (2006), pp.~221--252.

\bibitem{K2019}
{\sc E.~Kamenica}, {\em Bayesian persuasion and information design}, Annual
  Review of Economics, 11 (2019), pp.~249--272.

\bibitem{KEY2020}
{\sc X.-I. Kartala, N.~Englezos, and A.~N. Yannacopoulos}, {\em Future
  expectations modeling, random coefficient forward--backward stochastic
  differential equations, and stochastic viscosity solutions}, Mathematics of
  Operations Research, 45 (2020), pp.~403--433.

\bibitem{LL2007}
{\sc J.-M. Lasry and P.-L. Lions}, {\em Mean field games}, Japanese Journal of
  Mathematics, 2 (2007), pp.~229--260.

\bibitem{L2004}
{\sc A.~Lazrak}, {\em Generalized stochastic differential utility and
  preference for information}, The Annals of Applied Probability, 14 (2004),
  pp.~2149--2175.

\bibitem{LQ2003}
{\sc A.~Lazrak and M.~C. Quenez}, {\em A generalized stochastic differential
  utility}, Mathematics of operations research, 28 (2003), pp.~154--180.

\bibitem{MWZZ2015}
{\sc J.~MA, Z.~WU, D.~ZHANG, and J.~ZHANG}, {\em On well-posedness of
  forward-backward {SDEs}-a unified approach}, Annals of Applied Probability,
  25 (2015), pp.~2168--2214.

\bibitem{MY1999}
{\sc J.~Ma and J.~Yong}, {\em Forward-Backward Stochastic Differential
  Equations and Their Applications}, no.~1702, Springer Science \& Business
  Media, 1999.

\bibitem{MH2020}
{\sc Y.~Ma and M.~Huang}, {\em Linear quadratic mean field games with a major
  player: The multi-scale approach}, Automatica, 113 (2020), p.~108774.

\bibitem{MW1995}
{\sc M.~Miller and P.~Weller}, {\em Stochastic saddlepoint systems
  stabilization policy and the stock market}, Journal of Economic Dynamics and
  Control, 19 (1995), pp.~279--302.

\bibitem{NC2013}
{\sc M.~Nourian and P.~E. Caines}, {\em $\varepsilon$-nash mean field game
  theory for nonlinear stochastic dynamical systems with major and minor
  agents}, SIAM Journal on Control and Optimization, 51 (2013), pp.~3302--3331.

\bibitem{RR2006}
{\sc S.~T. Rachev and L.~R{\"u}schendorf}, {\em Mass Transportation Problems:
  Applications}, Springer Science \& Business Media, 2006.

\end{thebibliography}

\end{document}